\documentclass[11pt, reqno,a4paper]{article}
\usepackage{amsfonts,anysize,hyperref,amsmath,indentfirst,geometry}
\usepackage{bm}
\usepackage{amssymb}
\usepackage{mathrsfs}
\usepackage{amssymb}
\marginsize{25mm}{25mm}{15mm}{15mm}

\allowdisplaybreaks
\numberwithin{equation}{section}
\usepackage{graphicx}
\usepackage{amsmath}
\usepackage{amssymb}
\usepackage{amsfonts,amssymb,amsbsy,amsmath}
\usepackage{eqnarray}
\usepackage{diagbox}
\usepackage{mathptmx}
\usepackage{algorithm}
\usepackage{algorithmic}
\newtheorem{theorem}{Theorem}[section]
\newtheorem{definition}[theorem]{Definition}
\newtheorem{corollary}[theorem]{Corollary}

\newtheorem{remark}[theorem]{Remark}
\newtheorem{example}[theorem]{Example}
\allowdisplaybreaks

\begin{document}

\title{\Large \bfseries On the partial condition numbers for the indefinite least squares problem\thanks{The work is supported by the National Natural Science Foundation of China (Grant No. 11671060).}}
\author{{Hanyu Li\thanks{Corresponding author. Email addresses: lihy.hy@gmail.com or hyli@cqu.edu.cn; shaoxin.w@gmail.com}, Shaoxin Wang}}
\date{\small {College of Mathematics and Statistics, Chongqing University,
Chongqing, 401331, P. R. China \\
}} \maketitle
{\raggedleft\bfseries{\normalsize Abstract}}\\

The condition number of a linear function of the indefinite least squares solution is called the partial condition number for the indefinite least squares problem. In this paper, based on a new and very general condition number which can be called the unified condition number, we first present an expression of the partial unified condition number when the data space is measured by a general weighted product norm. Then, by setting the specific norms and weight parameters, we obtain the expressions of the partial normwise, mixed and componentwise condition numbers. Moreover, the corresponding structured partial condition numbers are also taken into consideration when the problem is structured. 
Considering the connections between the indefinite and total least squares problems, we derive the (structured) partial condition numbers for the latter, which generalize the ones in the literature.  To estimate these condition numbers effectively and reliably, the probabilistic spectral norm estimator and the small-sample statistical condition estimation method are applied and three related algorithms are devised. Finally, the obtained results are illustrated by numerical experiments.  \vspace{.3cm}

{\raggedleft \em AMS classification:}\ 65F20, 65F35, 65F30, 15A12, 15A60\\

{\raggedleft \em Keywords:} Indefinite least squares problem; Total least squares problem; Partial condition number; Normwise condition number; Mixed and componentwise condition number; Probabilistic spectral norm estimator; Small-sample statistical condition estimation

\section{Introduction}\label{sec.1}
\vspace{-2pt}
The indefinite least squares (ILS) problem is a generalization of the famous linear least squares (LLS) problem. It can be stated as follows:
\begin{equation}\label{1.1}
{\rm ILS:}\quad\mathop {\min }\limits_{x \in {\mathbb{R}^n}} {(b - Ax)^T}J(b - Ax),
\end{equation}
where $A \in {\mathbb{R}^{m \times n}}$ with $m\geq n$, $b \in \mathbb{R}^{m}$, and $J$ is a signature matrix defined as
\begin{equation*}\label{1.2}
J = \left[ {\begin{array}{*{20}{c}}
{{I_p} }&0\\
0&{- {I_q}}
\end{array}} \right], \quad p+q=m.
\end{equation*}
Hereafter, for any matrix $B$, $B^T$ denotes its transpose, and $\mathbb{R}^{n}$, $\mathbb{R}^{m \times n}$, and $I_r$ stand for the real vector space of dimension $n$, the set of $m\times n$ real matrices, and the identity matrix of order $r$, respectively. From  \cite{Boa,Cha}, it follows that the ILS problem \eqref{1.1} has a unique solution:
\[x(A,b) =M^{-1}{A^T}Jb \quad\ {\rm with }\ M=A^TJA\]
if and only if $A^{T}JA$ is positive definite. We will assume throughout this paper that the condition holds. Note that this condition implies that $p\geq n$ and $A(1:p,1:n)$  has full column rank and so does $A$ \cite{Boa}. So, for a genuinely ILS problem, $m>n$ is required.

The ILS problem was first proposed by Chandrasekaran et al. \cite{Cha} and finds many important applications in some areas. For example, it can be used to solve the total least squares (TLS) problem \cite{Huf}.  Also, we will encounter this problem in the area of optimization known as $H^{\infty}$-smoothing \cite{Has,Say}. The reader can refer to \cite{Cha} for the detailed explanations. So, some authors investigated its numerical algorithms, stability of algorithms, and perturbation analysis (e.g.,\cite{Boa,Cha,Grca,Liu11,Liu14,Mast,Wang,Xu11}). Considering that the condition number `plays a leading role in the study of both accuracy and complexity of numerical algorithms' \cite[p. vii]{Burg},  Bojanczyk et al.  \cite{Boa} and Grcar \cite{Grca} studied the normwise condition number of the ILS problem and presented an upper bound; Li et al.  \cite{Li14} discussed the mixed and componentwise condition numbers of this problem, and derived their explicit expressions and the easily computable upper bounds.

In this paper, by defining a unified condition number which includes the normwise, mixed and componentwise condition numbers as special cases, we mainly consider the partial condition numbers for the ILS problem when the data space $\mathbb{R}^{m \times n}\times \mathbb{R}^{m}$ is measured by a general weighted product norm. 

As mentioned in Abstract, the partial condition number is referred to the condition number of a linear function of the indefinite least squares solution $x(A,b)$, i.e., $L^Tx(A,b)$ with $L\in\mathbb{R}^{n \times k}$ ($k\leq n$). This kind of condition number was first studied by Cao and Petzold for linear systems based on the regular normwise condition number \cite{Cao}. Later, it was proposed for the LLS problem based on the normwise, mixed and componentwise condition numbers \cite{Ari,Bab09} and the TLS problem based on the normwise condition number \cite{Bab}. In \cite{Ari,Bab09,Bab,Cao}, the authors also provided some motivations for investigating this kind of condition number. For example, in practice, we may only be interested in the sensitivity of part of the elements of the solution and hence we only need to know the condition number of this part of the elements. The regular condition number cannot work well in this case. In addition, the regular condition number cannot evaluate the differences between the sensitivity of each element of the solution either. All of these problems can be tackled by the partial condition number since we can get the desired results by choosing different $L$. For example, when $L$ is the identity matrix or a column vector of the identity matrix, the partial condition number will reduce to the condition number of the solution $x(A,b)$ or of an element of the solution.

The general weighted product norm used to measure the data space $\mathbb{R}^{m \times n}\times \mathbb{R}^{m}$ in this paper is a generalization of the following weighted product norm
\begin{equation}\label{1.1000}
\left\| (\alpha A, \beta b) \right\|_F=\sqrt{\alpha^2\left\| A\right\|_F^2+\beta^2\left\| b \right\|_2},\quad \alpha>0, \beta>0,
\end{equation}
which was first used by Gratton for deriving the normwise condition number for the LLS problem \cite{Grat}. In \eqref{1.1000}, $\left\|\circ \right\|_F$ denotes the Frobenius norm of a matrix, and $\left\|\circ \right\|_2$ denotes the spectral norm of a matrix or the Euclidean norm of a vector. We will call the latter 2-norm uniformly later in this paper. Subsequently, the weighted product norm \eqref{1.1000} was applied to the partial normwise condition number for the LLS problem \cite{Ari} and the normwise condition number of the truncated singular value solution of a linear ill-posed problem \cite{Ber}. As pointed out in \cite{Grat}, this norm is very flexible. With it, we can monitor the perturbations on $A$ and $b$. For example, if $\alpha\rightarrow \infty$, no perturbation on $A$ will be permitted; similarly, if $\beta\rightarrow \infty$, there will be no perturbation on $b$ allowed. The norm \eqref{1.1000} was ever generalized to $\left\| (T A, \beta b) \right\|_F$  by Wei et al.  \cite{Wei} for studying the normwise condition number of the rank deficient LLS problem. Here, $T$ is a positive diagonal matrix. Later, the generalized norm was applied to the weighted LLS problem \cite{Yang}. The general weighted product norm of this paper is also a generalization of the above generalized norm; see the explanation following \eqref{2.2222}. So, in comparison, this kind of product norm has more advantages. 

Recently, the structured condition numbers of some problems such as the linear systems, the LLS problem, and the TLS problem have received a lot of attention. Rump \cite{Rump1,Rump2} presented the structured condition numbers of the linear systems with respect to normwise or componentwise distances. The obtained results generalized the corresponding ones in \cite{Hig}. Xu et al.  \cite{Xu} considered the structured normwise condition numbers for the LLS problem, while Cucker and Diao  \cite{Cucker} obtained its structured mixed and componentwise condition numbers. For the TLS problem, Li and Jia \cite{Lijia} derived its structured mormwise and mixed condition numbers. The results in \cite{Lijia,Rump1,Rump2,Xu} show that the structured condition number can be much tighter than the unstructured one in some cases. Like the structured condition numbers for the above problems, the structured partial condition numbers of the ILS problem are also of interest. We will investigate them in the fourth part of this paper corresponding to the results on the nonstructured partial condition numbers.

As introduced above or in \cite{Boa,Cha}, the ILS problem has a close relationship with the TLS problem. In fact, the TLS solution can be regarded as a solution to a special ILS problem; see \cite{Cha} or Section \ref{sec.5} below for details. In recent years, some authors studied the condition numbers of the TLS problem. Zhou et al. \cite{Zhou} considered the normwise, mixed, and componentwise condition numbers of the so called scaled TLS problem, a generalization of the TLS problem.
Afterward, Baboulin and Gratton \cite{Bab} investigated the partial normwise condition number of the TLS problem and provided some computable expressions. At the same time, Li and Jia \cite{Lijia} also presented an expression of the normwise condition number. The latest formula, and the lower and upper bounds of the normwise condition number for the TLS problem were given in \cite{Jia}. In addition, Xie et al. \cite{Xie14} showed that the three normwise condition numbers given in \cite{Bab, Lijia,Zhou} are mathematically equivalent. In the fifth part of this paper, we will find that the (structured) partial condition numbers of the TLS problem can be derived from the results of a special ILS problem. To our best knowledge, it is the first time to study the condition numbers for the TLS problem from the view of the ILS problem.

The rest of this paper is organized as follows. Section \ref{sec.2} presents some preliminaries. In Section \ref{sec.3}, we obtain the expressions of the partial unified condition number and the partial normwise, mixed and componentwise condition numbers of the ILS problem.
As mentioned above, Sections \ref{sec.4} and \ref{sec.5} are mainly devoted to the structured partial condition numbers of the ILS problem and the connections between the partial condition numbers of the ILS and TLS problems, respectively. Considering that computing a condition number may be expensive and a good estimate is acceptable for practical purpose \cite[Chapter 15]{Hig02}, in Section \ref{sec.6}, we provide the statistical estimates of the results derived in Sections \ref{sec.3} and \ref{sec.4} on basis of the probabilistic spectral norm estimator \cite{Hochs13} and the small-sample statistical condition estimation (SSCE) method \cite{Kenney94}. The numerical experiments for illustrating the obtained results are given in Section \ref{sec.7}. Finally, we present the conclusion of the whole paper.
\vspace{-6pt}

\section{Preliminaries}
\label{sec.2}
\vspace{-2pt}
Following \cite{Xieli}, we define the entry-wise division between the vectors $a\in \mathbb{R}^{p}$ and $b=[b_1,\cdots,b_p]^T\in \mathbb{R}^{p}$ by
\begin{eqnarray}
\label{eq.cdivis}
\frac{a}{b}= {\rm diag}^{\ddag}(b)a,
\end{eqnarray}
where ${\rm diag}^{\ddag}(b)$ is diagonal with diagonal elements $b_1^\ddag,\cdots,b_p^\ddag$. Here, for a number $c\in \mathbb{R}$, $c^{\ddag}$ is defined by
\begin{eqnarray*}
c^{\ddag}=\left\{\begin{array}{ll}
\frac{1}{c}, & \hbox{$c\neq 0$,} \\
1, & \hbox{$c= 0$ .}
\end{array}
\right.
\end{eqnarray*}
By \eqref{eq.cdivis}, we now define a new and general condition number.
\begin{definition}
\label{def.cond}
Let $F: \mathbb{R}^{p}\rightarrow \mathbb{R}^{q}$ be a continuous mapping defined on an open set $Dom(F)\in\mathbb{R}^{p}$, the domain of definition of $F$. Then the condition number of $F$ at $x\in Dom(F)$ is defined by
\begin{eqnarray*}
 \kappa_F(x)&=&\lim_{\delta\rightarrow 0}\sup_{0<\left\|\frac{\Delta x}{\beta}\right\|_\mu\leq \delta}\frac{\left\|\frac{F(x+\Delta x)-F(x)}{\xi}\right\|_{\nu}}{\left\|\frac{\Delta x}{\beta}\right\|_{\mu}},
\end{eqnarray*}
where $\|\cdot\|_\mu$ and $\|\cdot\|_\nu$ are the vector norms defined on $\mathbb{R}^{p}$ and $\mathbb{R}^{q}$, respectively, and $\beta\in \mathbb{R}^{p}$ and $\xi\in \mathbb{R}^{q}$ are parameters with a requirement that if some element of $\beta$ is zero, then the corresponding element of $\Delta x$ must be zero.
\end{definition}

\begin{remark}\label{rem1.gcond}{\rm
When we set $\beta$ to be the data $x$, the requirement on $\beta$ in Definition \ref{def.cond} means that the zero elements of $x$ do not perturb. As we know, in the floating point number system, a real number $\alpha$ can be represented as $fl(\alpha)=\alpha(1+\delta)$ with $|\delta|< \mu_0$, where $\mu_0$ is the unit roundoff \cite[p. 38]{Hig02}. Thus, when $\alpha=0$, we have $fl(\alpha)=0$. This fact shows that the zero element should not be perturbed and hence the mentioned requirement in Definition \ref{def.cond} is reasonable and acceptable.  }
\end{remark}
\begin{remark}\label{rem.gcond}{\rm
The condition number in Definition \ref{def.cond} can be called the unified condition number since it is very general
and covers several popular condition numbers. For example, when $\mu=\nu=2$, and $\beta=[\|x\|_2,\cdots,\|x\|_2]^T\in\mathbb{R}^{p}$ with $x\neq0$ and $\xi=[\|F(x)\|_2, \cdots, \|F(x)\|_2]^T\in\mathbb{R}^{q}$ with $F(x)\neq0$, we get the normwise condition number in \cite{Geu,Rice}; when $\mu=\nu=\infty$, and $\beta=x\neq0$ and $\xi=[\|F(x)\|_{\infty}, \cdots, \|F(x)\|_{\infty}]^T\in\mathbb{R}^{q}$ ($\xi=F(x)$) with $F(x)\neq0$, the mixed (componentwise) condition number in \cite{Gohb93,Xieli} follows. Moreover, the parameters $\beta$ and $\xi$ can be positive real numbers instead of vectors in Definition \ref{def.cond}. In this case, the entry-wise division between vectors reduces to the regular scalar multiplication between a scalar and a vector.
}\end{remark}

The operator `vec' and Kronecker product play important roles in obtaining the expression of the condition number. We introduce some necessary results on these two tools as follows.

For a matrix $A=[a_1,\cdots, a_n]\in
\mathbb{R}^{m\times n}$ with $a_i\in \mathbb{R}^{m}$, the operator 'vec' is defined as
\begin{equation*}\label{1.3}
{\rm vec}(A)=[a_1^T,\cdots,a_n^T]^T\in \mathbb{R}^{mn},
\end{equation*}
and the {Kronecker product} between $A =(a_{ij})\in {\mathbb{R}^{m \times n}}$ and $B \in
{\mathbb{R}^{p \times q}}$ is defined
by (e.g., \cite[Chapter 4]{Horn91}),
\[A \otimes B = \left[ {\begin{array}{*{20}{c}}
{{a_{11}}B}&{{a_{12}}B}& \cdots &{{a_{1n}}B}\\
{{a_{21}}B}&{{a_{22}}B}& \cdots &{{a_{2n}}B}\\
 \vdots & \vdots & \ddots & \vdots \\
{{a_{m1}}B}&{{a_{m2}}B}& \cdots &{{a_{mn}}B}
\end{array}} \right]\in {\mathbb{R}^{mp \times nq}}.\]
From the above definition, it is easy to find that when $m=1$ and $q=1$, i.e., when $A$ is a row vector and $B$ is a column vector,
\begin{eqnarray}\label{1.4}
A \otimes B=BA.
\end{eqnarray}
The following results on the operator `vec' and Kronecker product are from \cite[Chapter 4]{Horn91},
\begin{eqnarray}
&&(A \otimes B)^T= (A^T \otimes B^T),\label{1.5}\\
&&{\rm vec}(AXB) = \left({B^T} \otimes A\right){\rm vec}(X),\label{1.6}\\
&&\Pi_{mn} {\rm vec}(A) = {\rm vec}({A^T}),\label{1.7}\\
&&\Pi_{pm} (A \otimes B) \Pi_{nq}= (B \otimes A),\nonumber
\end{eqnarray}
where $X\in {\mathbb{R}^{n \times p}}$, and $\Pi_{st} \in {\mathbb{R}^{st \times st}}$ is the {\em vec-permutation matrix} which depends only on the dimensions $s$ and $t$. Note that if $n=1$, then $ \Pi_{nq}=I_q$ and hence
\begin{eqnarray}\label{1.8}
&&\Pi_{pm} (A \otimes B)= (B \otimes A).
\end{eqnarray}
In addition, from \cite[Chapter 4]{Horn91}, we also have
\begin{eqnarray}\label{1.9}
(A \otimes B)(C \otimes D)=(AC) \otimes (BD),
\end{eqnarray}
where the matrices $C$ and $D$ are of suitable orders.

\section{The partial condition numbers of the ILS problem}
\label{sec.3}
Let $L\in {\mathbb{R}^{n \times k}}$ with $k\leq n$ be a given matrix and be not perturbed numerically. We consider the following mapping
\begin{eqnarray*}
 g:\mathbb{R}^{m \times n}  \times \mathbb{R}^m &\rightarrow&\mathbb{R}^k\\
   (A,b)&\rightarrow&g(A,b) = L^T x(A,b) = L^T M^{ - 1} A^T Jb.
\end{eqnarray*}
From the discussions in \cite{Li14}, it follows that the mapping $g$ is continuously Fr\'{e}chet differentiable in a neighborhood of $(A,b)$. Denote by $g{'}(A,b)$ the Fr\'{e}chet derivative of $g$ at $(A,b)$. Thus, using the chain rules of composition of derivatives or from \cite{Boa,Li14}, we have
\begin{eqnarray}
 g{'}(A,b):\mathbb{R}^{m \times n}  \times \mathbb{R}^m  &\rightarrow& \mathbb{R}^k \nonumber \\
 (\Delta A,\Delta b) &\rightarrow& g{'} (A,b){\circ}{(\Delta A,\Delta b)} = L^T M^{ - 1} (\Delta A)^T Jr - L^T M^{ - 1} A^T J(\Delta A)x\nonumber \\
 &\quad& \quad\quad\quad\quad\quad\quad\quad\quad\quad + L^T M^{ - 1} A^T J(\Delta b), \label{2.1}
\end{eqnarray}
where $r=b-Ax$ and $g{'} (A,b){\circ}{(\Delta A,\Delta b)}$ denotes that we apply the mapping $g{'}(A,b)$ to the small perturbation variable $(\Delta A,\Delta b)$.
Then according to Definition \ref{def.cond} and the results in \cite{Geu, Rice}, and using the operator `vec', the condition number of $g$ at the point $(A,b)$ can be given by
\begin{eqnarray}\label{2.2222}
\kappa_{ILS}(A,b)
&=&\sup_{\left\|\mathrm{vec}\left(\frac{\Delta A}{\Psi},\frac{\Delta b}{\beta}\right)\right\|_{\mu}\neq 0}\frac{\left\|\frac{g' (A,b)\circ(\Delta A,\Delta b)}{\xi}\right\|_{\nu}}{\left\|\mathrm{vec}\left(\frac{\Delta A}{\Psi},\frac{\Delta b}{\beta}\right)\right\|_{\mu}},
\end{eqnarray}
where $\Psi\in\mathbb{R}^{m \times n}$, $\beta\in\mathbb{R}^{m}$ and $\xi\in\mathbb{R}^{k}$ are parameters with a requirement that if some element of $\Psi$ or $\beta$ is zero, then the corresponding element of $\Delta A$ or $\Delta b$ must be zero. As mentioned in Remark \ref{rem.gcond}, the parameters $\Psi$ and $\beta$ can be chosen to be positive real numbers. In this case, if we set $\mu=2$ further, the norm on the data space $\mathbb{R}^{m \times n}\times \mathbb{R}^{m}$ used in \eqref{2.2222} will reduce to the weighted product norm \eqref{1.1000}. What's more, if we set $\mu=2$, $\Psi$ to be a special positive matrix, and $\beta$ to be a positive real number, then the weighted product norm used in \cite{Wei} can be recovered. Consequently, the weighted product norm considered here is more general and hence has more advantages.

From the explanations in Section \ref{sec.1} and Remark \ref{rem.gcond}, we call the condition number $\kappa_{ILS} (A,b)$ the partial unified condition number of the ILS problem \eqref{1.1} with respect to $L$. An explicit expression of this condition number is presented as follows.
\begin{theorem}
\label{thm.gcils}
The partial unified condition number of the ILS problem \eqref{1.1} with respect to $L$ is
\begin{eqnarray}\label{2.2}
\kappa_{ILS} (A,b) =\left\|{\mathrm{diag}^\ddag(\xi)} {M_{g' } \mathrm{diag} (\mathrm{vec}(\Psi,\beta)) } \right\|_{\mu,\nu},
\end{eqnarray}
where
\begin{eqnarray}\label{2.3}
M_{g' }  = \left[ {{\left( {(Jr)^T  \otimes (L^T M^{ - 1} )} \right)\Pi _{mn}  - x^T  \otimes (L^T M^{ - 1} A^T J)},{L^T M^{ - 1} A^T J}} \right]
\end{eqnarray}
and $\|\cdot\|_{\mu,\nu}$ is the matrix norm induced by the vector norms $\|\cdot\|_\mu$ and $\|\cdot\|_\nu$.
\end{theorem}
\emph{Proof}. 
Applying the operator vec to $g{'} (A,b){\circ}{(\Delta A,\Delta b)}$ and using \eqref{1.6} and \eqref{1.7} gives
\begin{eqnarray}
 g' (A,b){\circ}(\Delta A,\Delta b) &=& {\rm vec}(g' (A,b){\circ}(\Delta A,\Delta b))\nonumber \\
  &=& \left( {(Jr)^T  \otimes (L^T M^{ - 1} )} \right)\Pi _{mn} {\rm vec}(\Delta A) \nonumber\\
 &\quad -& \left( {x^T  \otimes (L^T M^{ - 1} A^T J)} \right){\rm vec}(\Delta A)
  + L^T M^{ - 1} A^T J(\Delta b) \nonumber\\
&=& M_{g' }\left[ {\begin{array}{*{20}c}
   { {\rm vec}(\Delta A)}  \\
   { \Delta b}  \\
\end{array}} \right].\label{2.44}
\end{eqnarray}
Considering the requirement on $\Psi$ and $\beta$ in \eqref{2.2222}, we have
\begin{eqnarray}
\left[ {\begin{array}{*{20}c}
   { {\rm vec}(\Delta A)}  \\
   { \Delta b}  \\
\end{array}} \right]
  &=& \mathrm{diag} (\mathrm{vec}(\Psi,\beta))\left[ {\begin{array}{*{20}c}
   { {\rm vec}(\frac{\Delta A}{\Psi})}  \\
   { \frac{\Delta b}{\beta}}  \\
\end{array}} \right]. \label{2.4}
\end{eqnarray}
Substituting \eqref{2.4} into \eqref{2.44} and then into \eqref{2.2222} implies
\begin{eqnarray*}
\kappa_{ILS}(A,b)&=&\sup_{\left\|\mathrm{vec}\left(\frac{\Delta A}{\Psi},\frac{\Delta b}{\beta}\right)\right\|_{\mu}\neq 0}\frac{\left\|{\mathrm{diag}^\ddag(\xi)}M_{g' } \mathrm{diag} (\mathrm{vec}(\Psi,\beta))\left[ {\begin{array}{*{20}c}
  { {\rm vec}(\frac{\Delta A}{\Psi})}  \\
  { \frac{\Delta b}{\beta}}  \\
\end{array}} \right]\right\|_\nu}{\left\|\mathrm{vec}\left(\frac{\Delta A}{\Psi},\frac{\Delta b}{\beta}\right)\right\|_{\mu}}\nonumber\\
&=& \left\|{\mathrm{diag}^\ddag(\xi)}{M_{g' } \mathrm{diag} (\mathrm{vec}(\Psi,\beta)) } \right\|_{\mu,\nu}.\   \square
\end{eqnarray*}

Note that the expression of $\kappa_{ILS} (A,b)$ given in Theorem \ref{thm.gcils} is very general. In the following, we mainly concentrate on some specific norms and parameters to simplify and specify the expression.
\begin{theorem}[\bf 2-norm]
\label{thm.closedform}
When $\mu=\nu=2$, and the parameters $\Psi$, $\beta$, and $\xi$ are positive real numbers, the partial condition number \eqref{2.2} has the following two equivalent expressions
\begin{equation}\label{2.5}
\kappa_{2ILS} (A,b)  =  {\scalebox{1.3}{$\frac{\left\|L^T M^{ - 1} \left( {{\Psi^2\left\| r \right\|_2^2 }I_n  + \left( {\Psi^2{\left\| x \right\|_2^2 } + {{\beta ^2 }}} \right)A^T A - \Psi^2(xr^T A + A^T rx^T )} \right)M^{ - 1} L\right\|^{1/2}_2}{\xi}$}}
\end{equation}
and
\begin{equation}
\label{2.5'}
  \kappa_{2ILS} (A,b)=\frac{\left\|L^TM^{-1}\begin{bmatrix}
                         \Psi\|r\|_2(I_n-\frac{1}{\|r\|^2_2}A^Trx^T), & -\beta A^T, & \Psi \|x\|_2A^T(I_m-\frac{1}{\|r\|_2^2}rr^T) \\
                       \end{bmatrix}\right\|_2}{\xi}.
\end{equation}
\end{theorem}
\emph{Proof}.  Under the hypothesis of this theorem, from Theorem \ref{thm.gcils}, we have
\begin{eqnarray}
\kappa_{2ILS} (A,b)=\frac{\left\|\left[\Psi M_1,\beta {L^T M^{ - 1} A^T J}\right] \right\|_2 }{\xi} ,\label{2.66}
\end{eqnarray}
where $M_1  = \left( {(Jr)^T  \otimes (L^T M^{ - 1} )} \right)\Pi _{mn}  - x^T  \otimes (L^T M^{ - 1} A^T J).$
Note that, for any matrix $X\in \mathbb{R}^{m\times n}$, $\left\| X \right\|_2  = \left\| X X^T  \right\|_2^{1/2}$. Thus,
\begin{eqnarray}
\kappa_{2ILS} (A,b)&=& \frac{\left\|\Psi^2{M_1 M_1^T } + \beta^2 {L^T M^{ - 1} A^T AM^{ - 1} L} \right\|_2^{1/2}}{\xi}.\label{2.6}
\end{eqnarray}
Considering \eqref{1.5}, \eqref{1.9}, \eqref{1.8}, and \eqref{1.4}, we obtain
\begin{align*}
&M_1 M_1^T = \left( {\left( {(Jr)^T  \otimes (L^T M^{ - 1} )} \right)\Pi _{mn}  - x^T  \otimes (L^T M^{ - 1} A^T J)} \right)\\
&\quad\quad\quad\   \times  \left( {\Pi _{mn}^T \left( {(Jr) \otimes (M^{ - 1} L)} \right)
  - x \otimes (JAM^{ - 1} L)} \right)\quad \textrm{ by \eqref{1.5}} \\
 &\quad\quad\  \  \ = \left( {\left( {(Jr)^T (Jr)} \right) \otimes \left( {L^T M^{ - 2} L} \right)} \right) + \left( {(x^T x) \otimes (L^T M^{ - 1} A^T AM^{ - 1} L)} \right)\\
 &\quad\quad\quad\   - ( {(Jr)^T  \otimes (L^T M^{ - 1} )} )( {(JAM^{ - 1} L) \otimes x} ) - ( {x^T  \otimes (L^T M^{ - 1} A^T J)} )( {(M^{ - 1} L) \otimes (Jr)} )   \textrm{ by \eqref{1.9} and \eqref{1.8}}\\
 &\quad\quad\  \  \ = \left( {\left( {(Jr)^T (Jr)} \right) \otimes \left( {L^T M^{ - 2} L} \right)} \right) + \left( {(x^T x) \otimes (L^T M^{ - 1} A^T AM^{ - 1} L)} \right)\\
 &\quad\quad\quad\   - \left( {(r^T AM^{ - 1} L) \otimes (L^T M^{ - 1} x)} \right) - \left( {(x^T M^{ - 1} L) \otimes (L^T M^{ - 1} A^T r)} \right)  \quad \textrm{ by \eqref{1.9}}\\
 &\quad\quad\  \  \ = \left\| r \right\|_2^2  {L^T M^{ - 2} L}  + \left\| x \right\|_2^2 L^T M^{ - 1} A^T AM^{ - 1} L
  - L^T M^{ - 1} xr^T AM^{ - 1} L\\
  &\quad\quad\quad\   - L^T M^{ - 1} A^T rx^T M^{ - 1} L.  \quad \textrm{ by \eqref{1.4}}
\end{align*}
Substituting the above equality into \eqref{2.6} gives \eqref{2.5}.

On the other hand, if we set
$$K=\begin{bmatrix}
                         \Psi\|r\|_2(I_n-\frac{1}{\|r\|^2_2}A^Trx^T), & -\beta A^T, & \Psi \|x\|_2A^T(I_m-\frac{1}{\|r\|_2^2}rr^T) \\
                       \end{bmatrix},$$
we can check that
\begin{equation*}
   \Psi^2{M_1 M_1^T } + \beta^2 {L^T M^{ - 1} A^T AM^{ - 1} L} =L^TM^{-1}KK^TM^{-1}L.
\end{equation*}
Again, by the equality $\left\| X \right\|_2  = \left\| X X^T  \right\|_2^{1/2}$ and \eqref{2.6}, we have \eqref{2.5'}.
$\square$

\begin{remark}{\rm
Note that the orders of the matrices in \eqref{2.5}, \eqref{2.5'}, and \eqref{2.66} are $k\times k$, $k\times (2m+n)$, and $k\times (mn+m)$, respectively. Hence, when $m$ and $n$ are very large, both of the expressions \eqref{2.5} and \eqref{2.5'} reduce the storage requirements significantly. However, forming the matrix in \eqref{2.5} explicitly is not desirable because computing the cross product $A^TA$ may be potentially unstable \cite[p. 386]{Hig02}. Therefore, in comparison, the expression \eqref{2.5'} seems to be more preferred.
}\end{remark}

\begin{remark}{\rm
The condition number in Theorem \ref{thm.closedform} is the simplified partial normwise condition number of the ILS problem. Setting $L=I_n$ and $\Psi=\beta=\xi=1$ in \eqref{2.66}, and using the property on the spectral norm that for the matrices $C$ and $D$ of suitable orders, $\left\|[C,D]\right\|_2\leq\left\|C\right\|_2+\left\|D\right\|_2$, we have
\begin{eqnarray*}
\kappa_{2ILS} (A,b) &\leq&   \left\|{\left( {(Jr)^T  \otimes M^{ - 1} } \right)\Pi _{mn}  - x^T  \otimes (M^{ - 1} A^T J)}\right\|_2+ \left\|{M^{ - 1} A^T }\right\|_2,
\end{eqnarray*}
which is equivalent to the upper bound of the normwise condition number for the ILS problem given in \cite[(2.10)]{Boa} or  \cite[(4.5)]{Grca} in essence. 
}\end{remark}

As mentioned in Section \ref{sec.1}, the ILS problem is a generalization of the LLS problem. Thus, setting $J=I_{m}$ in the above results and noting, in this case, $M=A^TA$ and $A^Tr=0$, we have the corresponding results on the partial condition numbers for the LLS problem.

\begin{corollary}
The partial unified condition number of the LLS problem with respect to $L$ is
\begin{eqnarray}\label{2.2222222}
\kappa_{LLS} (A,b) =\left\|{\mathrm{diag}^\ddag(\xi)} {\widetilde{M}_{\widetilde{g}' } \mathrm{diag} (\mathrm{vec}(\Psi,\beta)) } \right\|_{\mu,\nu},
\end{eqnarray}
where
\begin{eqnarray}\label{2.333333}
\widetilde{M}_{\widetilde{g}' }  = \left[ {{\left( {r^T  \otimes (L^T (A^TA)^{ - 1} )} \right)\Pi _{mn}  - x^T  \otimes (L^T (A^TA)^{ - 1} A^T )},{L^T (A^TA)^{ - 1} A^T }} \right].
\end{eqnarray}

If $\mu=\nu=2$, and the parameters $\Psi$, $\beta$, and $\xi$ are positive real numbers, then
\begin{eqnarray}
\kappa_{2LLS} (A,b) &=& \frac{ \left\|{\Psi^2\left\| r \right\|_2^2 } L^T (A^TA)^{-2} L + \left( {{\Psi^2\left\| x \right\|_2^2 } + {{\beta ^2 }}} \right)L^T (A^TA)^{-1}L\right\|_2^{1/2}}{\xi} \label{2.8}
\end{eqnarray}
and
\begin{equation}
\label{2.8'}
  \kappa_{2LLS} (A,b)=\frac{\left\|L^T(A^TA)^{-1}\begin{bmatrix}
                         \Psi\|r\|_2I_n, & -\beta A^T, & \Psi \|x\|_2A^T \\
                       \end{bmatrix}\right\|_2}{\xi}.
\end{equation}
\end{corollary}

\begin{remark}{\rm 
If $L$ is a column vector, i.e., $k=1$, then \eqref{2.8} reduces
\begin{eqnarray*}
  \kappa_{2LLS} (A,b)  &=& \left({\Psi^2{\left\| r \right\|_2^2 }} \left\|L^T (A^TA)^{-1} \right\|_2^2 + \left( {{\Psi^2\left\| x \right\|_2^2 } + {{\beta ^2 }}} \right)\left\|L^T A^{\dag}\right\|_2^2\right)^{1/2},\label{2.9}
\end{eqnarray*}
which is just the result given in \cite[Corollary 1]{Ari}. Hereafter, $A^{\dag}$ denotes the Moore-Penrose inverse of the matrix $A$ (e.g., \cite{Hig02}).}
\end{remark}

\begin{remark}{\rm
Let $A = U\Sigma V^T$ be the thin singular value decomposition of $A$ appearing in the LLS problem with $U \in \mathbb{R}^{m \times n} ,V \in \mathbb{R}^{n \times n} $, and
$\Sigma  = {\rm diag}(\sigma _1 , \cdots ,\sigma _n )$ satisfying $U^T U = I_n  = V^T V = VV^T $ and $\sigma _1  \ge  \cdots  \ge \sigma _n  > 0$ (e.g., \cite{Hig02}).
Then $A^T A = V\Sigma ^2 V^T $. Substituting this equation into \eqref{2.8} and \eqref{2.8'} yields
\begin{equation}
\label{3.14}
\kappa_{2LLS} (A,b) = \frac{\left\| {L^T V\Sigma ^{ - 2} S^2 \Sigma ^{ - 2} V^T L} \right\|_2^{1/2}}{\xi}  = \frac{\left\| {S\Sigma ^{ - 2} V^T L} \right\|_2}{\xi}
\end{equation}
and
\begin{equation}
\label{3.15}
\kappa_{2LLS} (A,b)=\frac{\left\|L^TV\begin{bmatrix}
                         \Psi\|r\|_2\Sigma^{-2}, & -\beta \Sigma^{-1}, & \Psi \|x\|_2 \Sigma^{-1} \\
                       \end{bmatrix}P^T\right\|_2}{\xi},
\end{equation}
where $S$ is a diagonal matrix with the diagonal elements
$$S_{ii}  = \sqrt {\Psi^2{\left\| r \right\|_2^2 } + \left( {\Psi^2{\left\| x \right\|_2^2 } +{{\beta ^2 }}} \right)\sigma _i^2 } ,\quad i = 1, \cdots ,n,$$
and $P$ is a column orthnormal and block-diagonal matrix with $V$, $U$, and $U$ on its diagonal.

When $\xi=1$, \eqref{3.14} is just the expression given in \cite[Theorem 1]{Ari}, where it was derived by an alternative approach. Although it is easy to check that  \eqref{3.15} is equivalent to \eqref{3.14}, the expression \eqref{3.15} (or \eqref{2.8'}) is new as far as we know.

In addition, since $M=A_p^TA_p-A_q^TA_q$ if $A$ is divided into $
A = \left[A_p^T ,\ A_q^T \right]^T$ with $A_p \in \mathbb{R}^{p \times n}$ and $A_q \in \mathbb{R}^{q \times n}$, we can apply the generalized singular value decomposition of the matrix pair $A_p, A_q$ \cite{Paige} to rewrite the condition number \eqref{2.5} or \eqref{2.5'} of the ILS problem as done above for the  LLS problem. Due to its complexity and length, the topic will be considered in a separate paper.
}\end{remark}

Now we consider the partial condition number with $\mu=\nu=\infty$ for the ILS problem \eqref{1.1}, from which the partial mixed and componentwise condition numbers follow.
\begin{theorem}[\bf $\infty$-norm]
\label{thm.m&c}
When $\mu=\nu=\infty$, the partial condition number of the ILS problem \eqref{1.1} with respect to $L$ is
\begin{eqnarray}\label{eq.m&c}
\kappa_{{\infty}ILS} (A,b) =\left\|{{\rm diag}^\ddag(\xi)}{M_{g' } \mathrm{diag} (\mathrm{vec}(\Psi,\beta)) } \right\|_{\infty}=\left\| \frac{{\left|M_{g' }\right|  \left|\mathrm{vec}(\Psi,\beta)\right| }}{\left|\xi\right|} \right\|_{\infty},
\end{eqnarray}
where $M_{g' }$ is defined by \eqref{2.3}.

In particular, setting $\Psi=A$, $\beta=b$, and $\xi=\left[\|L^Tx(A,b)\|_{\infty},\cdots,\|L^Tx(A,b)\|_{\infty}\right]^T$ or $\xi=L^Tx(A,b)$, we get the corresponding partial mixed or componentwise condition number
\begin{eqnarray}
\label{eq.mcdils}
\kappa_{mILS} (A,b) =\frac{\left\|{\left|M_{g' }\right|  \left|\mathrm{vec}(A,b)\right| } \right\|_{\infty}}{\|L^Tx(A,b)\|_{\infty}}
\end{eqnarray}
or
\begin{eqnarray}
\label{eq.ccdils}
\kappa_{cILS} (A,b) =\left\| \frac{{\left|M_{g' }\right|  \left|\mathrm{vec}(A,b)\right| }}{|L^Tx(A,b)|} \right\|_{\infty}.
\end{eqnarray}
\end{theorem}
\emph{Proof}.  Letting $\mu=\nu=\infty$ in \eqref{2.2} gives the first part of \eqref{eq.m&c}. The second part of \eqref{eq.m&c} can be obtained by considering the proof of Lemma 2 in \cite{Cucker07} and \eqref{eq.cdivis}. The expressions \eqref{eq.mcdils} and  \eqref{eq.ccdils} are the straightforward results of \eqref{eq.m&c}.
$\square$

\begin{remark}{\rm
When $L=I_n$, \eqref{eq.mcdils} and  \eqref{eq.ccdils} reduce to the mixed and componentwise condition numbers for the ILS problem \eqref{1.1} established in \cite{Li14}. If we set $J=I_m$ in Theorem \ref{thm.m&c}, we have the corresponding results for the LLS problem:
\begin{eqnarray*}
\kappa_{{\infty}LLS} (A,b) =\left\|{{\rm diag}^\ddag(\xi)}{\widetilde{M}_{\widetilde{g}' } \mathrm{diag} (\mathrm{vec}(\Psi,\beta)) } \right\|_{\infty}=\left\| \frac{{\left|\widetilde{M}_{\widetilde{g}' }\right|  \left|\mathrm{vec}(\Psi,\beta)\right| }}{\left|\xi\right|} \right\|_{\infty},
\end{eqnarray*}
\begin{eqnarray*}
\kappa_{mLLS} (A,b) =\frac{\left\| \left| {\left( {r^T  \otimes (L^T (A^TA)^{ - 1} )} \right)\Pi _{mn}  - x^T  \otimes (L^T (A^TA)^{ - 1} A^T)}\right| \left|\mathrm{vec}(A)\right|+\left|{L^T (A^TA)^{ - 1} A^T} \right||b| \right\|_{\infty}}{\|L^Tx(A,b)\|_{\infty}},
\end{eqnarray*}
and
\begin{eqnarray*}
\kappa_{cLLS} (A,b) =\left\| \frac{{\left| {\left( {r^T  \otimes (L^T (A^TA)^{ - 1} )} \right)\Pi _{mn}  - x^T  \otimes (L^T (A^TA)^{ - 1} A^T)}\right| \left|\mathrm{vec}(A)\right|+\left|{L^T (A^TA)^{ - 1} A^T} \right||b| }}{|L^Tx(A,b)|} \right\|_{\infty},
\end{eqnarray*}
where $\widetilde{M}_{\widetilde{g}' }$ is defined by \eqref{2.333333}, and the mixed and componentwise condition numbers are the same as the ones derived from \cite{Bab09}.
}\end{remark}
\vspace{-6pt}

\section{The  structured partial condition numbers of the ILS problem}
\label{sec.4}
\vspace{-2pt}
Let $\mathbb{S}_1\subseteq \mathbb{R}^{m\times n}$ and $\mathbb{S}_2\subseteq \mathbb{R}^{m}$ be two linear subspaces. The former consists of a class of matrices having the same structure such as the symmetric matrices, the Toeplitz matrices, the Hankel matrices and so on (e.g., \cite{Cucker, Hig, Rump1, Rump2}), and the latter comprises a class of structured vectors. According to \cite{Hig,Lijia,Rump1}, we have that if $A\in \mathbb{S}_1$ and $b \in \mathbb{S}_2$, then
\begin{eqnarray*}
{\rm vec}(A)=\Phi_{\mathbb{S}_1}s_1,\quad b=\Phi_{\mathbb{S}_2}s_2,
\end{eqnarray*}
where $\Phi_{\mathbb{S}_1}\in \mathbb{R}^{mn\times k_1}$ and $\Phi_{\mathbb{S}_2}\in \mathbb{R}^{m\times k_2}$ are the fixed structure matrices respectively reflecting the structures of $\mathbb{S}_1$ and $\mathbb{S}_2$, and $s_1\in \mathbb{R}^{k_1}$ and $s_2\in \mathbb{R}^{k_2}$ are the vectors of independent parameters in the structured matrices and vectors, respectively.

Now, in a similar manner as \eqref{2.2222}, we present the definition of the structured partial unified condition number for the ILS problem \eqref{1.1} with respect to $L$:
\begin{eqnarray}
\label{def.scd}
\kappa^S_{ILS}(A,b)&=&\mathop {\sup }\limits_{{\left\|\mathrm{vec}\left(\frac{\Delta A}{\Psi},\frac{\Delta b}{\beta}\right)\right\|_{\mu}\neq 0}\hfill \atop
  \scriptstyle
\Delta A,\Psi\in \mathbb{S}_1, \Delta b,\beta\in \mathbb{S}_2\hfill}\frac{\left\|\frac{g' (A,b)\circ(\Delta A,\Delta b)}{\xi}\right\|_{\nu}}{\left\|\mathrm{vec}\left(\frac{\Delta A}{\Psi},\frac{\Delta b}{\beta}\right)\right\|_{\mu}},
\end{eqnarray}
where the requirement on the parameters $\Psi$ and $\beta$ is the same as the one in \eqref{2.2222}. Moreover, two additional requirements on $\Psi$ and $\beta$, i.e., $\Psi\in \mathbb{S}_1$ and $\beta\in\mathbb{S}_2$, are added. Since it is difficult to take supremum over the structured data space \cite{Rump1}, in the following, we mainly focus on some specific norms to tackle this problem.

Firstly, we set $\mu=\nu=2$. In this case, substituting \eqref{2.4} and \eqref{2.44} into \eqref{def.scd} gives
\begin{align}\label{4.3}
\kappa^S_{2ILS}(A,b)&=&\mathop {\sup }\limits_{{\left\|\mathrm{vec}\left(\frac{\Delta A}{\Psi},\frac{\Delta b}{\beta}\right)\right\|_{2}\neq 0}\hfill \atop
  \scriptstyle
\Delta A,\Psi\in \mathbb{S}_1, \Delta b,\beta\in \mathbb{S}_2\hfill}\frac{\left\|{\rm diag}^\ddag({\xi}){ M_{g' } \mathrm{diag} (\mathrm{vec}(\Psi,\beta))\left[ {\begin{array}{*{20}c}
   { {\rm vec}(\frac{\Delta A}{\Psi})}  \\
   { \frac{\Delta b}{\beta}}  \\
\end{array}} \right]}\right\|_{2}}{\left\|\mathrm{vec}\left(\frac{\Delta A}{\Psi},\frac{\Delta b}{\beta}\right)\right\|_{2}}.
\end{align}
Since $\Delta A,\Psi\in \mathbb{S}_1$ and $\Delta b,\beta\in \mathbb{S}_2$, based on the explanation at the beginning of this section, we have
\begin{eqnarray}\label{4.200000}
{\rm vec}(\Delta A)=\Phi_{\mathbb{S}_1}\Delta s_1,\quad\mathrm{vec}(\Psi)=\Phi_{\mathbb{S}_1}\varphi,\quad \Delta b=\Phi_{\mathbb{S}_2}\Delta s_2, \quad \beta=\mathrm{vec}(\beta)=\Phi_{\mathbb{S}_2}\theta,
\end{eqnarray}
where $\Delta s_1\in \mathbb{R}^{k_1}$ and $\Delta s_2\in \mathbb{R}^{k_2}$ can be regarded as the corresponding perturbations of $s_1$ and $s_2$, and $\varphi\in \mathbb{R}^{k_1}$ and $\theta\in \mathbb{R}^{k_2}$ can be interpreted as the vectors of the parameters that really work in $\Psi$ and $\beta$.
As a result,
\begin{eqnarray*}
{\rm vec}\left(\frac{\Delta A}{\Psi}\right)=\Phi_{\mathbb{S}_1}\frac{\Delta s_1}{\varphi},\quad \frac{\Delta b}{\beta}=\Phi_{\mathbb{S}_2}\frac{\Delta s_2}{\theta},
\end{eqnarray*}
which can be written together as
\begin{eqnarray}\label{4.400}
 \left[ {\begin{array}{*{20}c}
   { {\rm vec}(\frac{\Delta A}{\Psi})}  \\
   { \frac{\Delta b}{\beta}}  \\
\end{array}} \right]=\left[ {\begin{array}{*{20}c}
   \Phi_{\mathbb{S}_1} &0  \\
   0&\Phi_{\mathbb{S}_2}  \\
\end{array}} \right]\left[ {\begin{array}{*{20}c}
   \frac{\Delta s_1}{\varphi} \\
   \frac{\Delta s_2}{\theta} \\
\end{array}} \right].
\end{eqnarray}
Substituting the above equation into \eqref{4.3} yields
{\small\begin{equation}\label{4.4}
\kappa^S_{2ILS}(A,b)=\mathop {\sup }\limits_{\tiny{\left\|\left[ {\begin{array}{*{20}c}
   \Phi_{\mathbb{S}_1} &0  \\
   0&\Phi_{\mathbb{S}_2}  \\
\end{array}} \right]\left[ {\begin{array}{*{20}c}
   \frac{\Delta s _1}{\varphi} \\
  \frac{\Delta s_2}{\theta} \\
\end{array}} \right]\right\|_{2}\neq 0}\hfill \atop
\quad \Delta s_1,\varphi\in \mathbb{R}^{k_1},\Delta s_2,\theta\in \mathbb{R}^{k_2}\hfill}\frac{\left\|{\rm diag}^\ddag({\xi}){M_{g' } \mathrm{diag} (\mathrm{vec}(\Psi,\beta))\left[ {\begin{array}{*{20}c}
   \Phi_{\mathbb{S}_1} &0  \\
   0&\Phi_{\mathbb{S}_2}  \\
\end{array}} \right]\left[ {\begin{array}{*{20}c}
     \frac{\Delta s _1}{\varphi} \\
    \frac{\Delta s_2}{\theta} \\
\end{array}} \right]}\right\|_{2}}{\left\|\left[ {\begin{array}{*{20}c}
   \Phi_{\mathbb{S}_1} &0  \\
   0&\Phi_{\mathbb{S}_2}  \\
\end{array}} \right]\left[ {\begin{array}{*{20}c}
   \frac{\Delta s _1}{\varphi} \\
  \frac{\Delta s_2}{\theta} \\
\end{array}} \right]\right\|_{2}}.
\end{equation}}
 Since
\begin{eqnarray*}
{{\left\| {\left[ {\begin{array}{*{20}c}
   \Phi_{\mathbb{S}_1} &0  \\
   0&\Phi_{\mathbb{S}_2}  \\
\end{array}} \right]\left[ {\begin{array}{*{20}c}
   \frac{\Delta s _1}{\varphi} \\
  \frac{\Delta s_2}{\theta} \\
\end{array}} \right]} \right\|_2 }}=\left\|\left[ {\begin{array}{*{20}c}
   \frac{\Delta s _1}{\varphi} \\
  \frac{\Delta s_2}{\theta} \\
\end{array}} \right]^T {\left[ {\begin{array}{*{20}c}
   \Phi_{\mathbb{S}_1}^T\Phi_{\mathbb{S}_1} &0  \\
   0&\Phi_{\mathbb{S}_2}^T\Phi_{\mathbb{S}_2}  \\
\end{array}} \right]\left[ {\begin{array}{*{20}c}
   \frac{\Delta s _1}{\varphi} \\
  \frac{\Delta s_2}{\theta} \\
\end{array}} \right]} \right\|_2^{1/2}
\end{eqnarray*}
and the structured matrices $\Phi_{\mathbb{S}_1}$ and $\Phi_{\mathbb{S}_2}$ are column orthogonal \cite{Lijia},
\begin{eqnarray}\label{4.5}
{{\left\| {\left[ {\begin{array}{*{20}c}
   \Phi_{\mathbb{S}_1} &0  \\
   0&\Phi_{\mathbb{S}_2}  \\
\end{array}} \right]\left[ {\begin{array}{*{20}c}
   \frac{\Delta s _1}{\varphi} \\
  \frac{\Delta s_2}{\theta} \\
\end{array}} \right]} \right\|_2 }}=\left\| {\left[ {\begin{array}{*{20}c}
   D_1 &0  \\
   0&D_2  \\
\end{array}} \right]\left[ {\begin{array}{*{20}c}
   \frac{\Delta s _1}{\varphi} \\
  \frac{\Delta s_2}{\theta} \\
\end{array}} \right]} \right\|_2,
\end{eqnarray}
where $D_1={\rm diag}(w_1)$ and $D_2={\rm diag}(w_2)$ with
\begin{eqnarray*}
w_1=[\left\|\Phi_{\mathbb{S}_1}(:,1)\right\|_2,\cdots,\left\|\Phi_{\mathbb{S}_1}(:,k_1)\right\|_2]^T,\quad w_2=[\left\|\Phi_{\mathbb{S}_2}(:,1)\right\|_2,\cdots,\left\|\Phi_{\mathbb{S}_2}(:,k_2)\right\|_2]^T.
\end{eqnarray*}
Combining \eqref{4.4} and \eqref{4.5} implies
\begin{align*}
&{\scalebox{1}{$\kappa_{2ILS}^S (A,b)=\mathop {\sup }\limits_{\tiny{{\left\| {\left[ {\begin{array}{*{20}c}
   D_1 &0  \\
   0&D_2  \\
\end{array}} \right]\left[ {\begin{array}{*{20}c}
   \frac{\Delta s _1}{\varphi} \\
  \frac{\Delta s_2}{\theta} \\
\end{array}} \right]} \right\|_2 }\neq 0}\hfill \atop
\quad \Delta s_1,\varphi\in \mathbb{R}^{k_1},\Delta s_2,\theta\in \mathbb{R}^{k_2}\hfill}\frac{{\left\| {\rm diag}^\ddag({\xi}){{M_{g' }  \mathrm{diag} (\mathrm{vec}(\Psi,\beta))\left[ {\begin{array}{*{20}c}
   \Phi_{\mathbb{S}_1}D_1^{-1} &0  \\
   0&\Phi_{\mathbb{S}_2} D_2^{-1} \\
\end{array}} \right]\left[ {\begin{array}{*{20}c}
   D_1 &0  \\
   0&D_2  \\
\end{array}} \right]\left[ {\begin{array}{*{20}c}
   \frac{\Delta s _1}{\varphi} \\
  \frac{\Delta s_2}{\theta} \\
\end{array}} \right]}}\right\|_2 }}{{\left\| {\left[ {\begin{array}{*{20}c}
   D_1 &0  \\
   0&D_2  \\
\end{array}} \right]\left[ {\begin{array}{*{20}c}
   \frac{\Delta s _1}{\varphi} \\
  \frac{\Delta s_2}{\theta} \\
\end{array}} \right]} \right\|_2 }}$}} .
\end{align*}
Considering that $\left[ {\begin{array}{*{20}c}
   D_1 &0  \\
   0&D_2  \\
\end{array}} \right]$ is nonsingular, for the above equation, we can take the supremum over all the parameters in $\mathbb{R}^{k_1+k_2}$, 
and hence get the following theorem.

\begin{theorem}[\bf Structured 2-norm]
Let $A\in \mathbb{S}_1$ and $b \in \mathbb{S}_2$. Then the structured partial condition number under 2-norm of the ILS problem \eqref{1.1} with respect to $L$ is
\begin{eqnarray}\label{4.2222}
\kappa_{2ILS}^S (A,b) = \left\|{\rm diag}^\ddag({\xi}){M_{g' }  \mathrm{diag} (\mathrm{vec}(\Psi,\beta))\left[ {\begin{array}{*{20}c}
   \Phi_{\mathbb{S}_1}D_1^{-1} &0  \\
   0&\Phi_{\mathbb{S}_2} D_2^{-1} \\
\end{array}} \right]} \right\|_2,
\end{eqnarray}
where $M_{g' }$ is defined by \eqref{2.3}.

If the parameters $\Psi$, $\beta$, and $\xi$ are positive real numbers, then
\begin{align}
&\kappa_{2ILS}^S (A,b)\nonumber \\
&=\frac{\left\|\left[ {{\Psi\left(\left( {(Jr)^T  \otimes (L^T M^{ - 1} )} \right)\Pi _{mn}  - x^T  \otimes (L^T M^{ -1}A^TJ)\right)}\Phi_{\mathbb{S}_1}D_1^{-1} ,\beta{L^TM^{ -1}A^TJ }}\Phi_{\mathbb{S}_2}D_2^{-1}  \right]\right\|_2}{\xi}.\label{4.22222}
\end{align}
\end{theorem}

\begin{remark}\label{Remark}{\rm
It is easy to verify that
$$ \left[ {\begin{array}{*{20}c}
   \Phi_{\mathbb{S}_1}D_1^{-1} &0  \\
   0&\Phi_{\mathbb{S}_2} D_2^{-1} \\
\end{array}} \right]$$
is column orthonormal. Thus,
$$\kappa_{2ILS}^S (A,b)\leq \kappa_{2ILS} (A,b).$$
That is, the structured partial condition number \eqref{4.22222} is always tighter than the unstructured one \eqref{2.66}. This fact can also be seen from the definitions of the condition numbers \eqref{def.scd} and \eqref{2.2222}.  As done in \cite{Rump1, Rump2, Xu}, it is interesting to investigate the ratio $\kappa_{2ILS}^S (A,b)/\kappa_{2ILS} (A,b)$ to see whether the former can be much smaller than the latter. We won't consider this topic in this paper, and only provide a numerical example in
Section \ref{sec.7} to show that the structured partial condition numbers, including the structured partial mixed and componentwise condition numbers given below, are indeed tighter than the corresponding unstructured ones.
}\end{remark}

\begin{remark}{\rm
Setting $J=I_m$ in \eqref{4.22222}, we have the structured partial condition number under 2-norm of the LLS problem with respect to $L$:
\begin{align}
\kappa_{2LLS}^S (A,b) = \frac{\left\|\left[ {{\Psi\left(\left( {r^T  \otimes (L^T (A^TA)^{ - 1} )} \right)\Pi _{mn}  - x^T  \otimes (L^T A^{ \dag})\right)}\Phi_{\mathbb{S}_1}D_1^{-1} ,\beta{L^T A^{ \dag} }}\Phi_{\mathbb{S}_2}D_2^{-1}  \right]\right\|_2}{\xi} .\label{4.7}
\end{align}
Further, the corresponding result for the linear system can be obtained by setting $A$ to be nonsingular and noting $r=0$:
\begin{eqnarray}
\kappa_{2LLS}^S (A,b) &=& \frac{\left\|\left[ {{ \Psi\left( -  x^T  \otimes (L^T A^{ -1})\right)}\Phi_{\mathbb{S}_1}D_1^{-1} ,\beta{L^T A^{ -1} }}\Phi_{\mathbb{S}_2}D_2^{-1}  \right]\right\|_2}{\xi}. \label{4.70}
\end{eqnarray}
Noting that
\begin{align*}
\frac{\left\| {{ \Psi\left( -  x^T  \otimes (L^T A^{ -1})\right)}}\Phi_{\mathbb{S}_1}\right\|_2}{\left\| D_1 \right\|_2}\leq\left\| {{\Psi\left(-  x^T  \otimes (L^T A^{ -1})\right)}}\Phi_{\mathbb{S}_1}D_1^{-1} \right\|_2\leq\left\| {{ \Psi\left( -  x^T  \otimes (L^T A^{ -1})\right)}}\Phi_{\mathbb{S}_1}\right\|_2\left\| D_1^{-1} \right\|_2,
\end{align*}
we can find that the structured condition number \eqref{4.70} with $L=I_n$ and $\Psi=\beta=\xi=1$ is equivalent to the one given in \cite{Rump1} in essence.

In addition, the structured condition number for the LLS problem under two conditions derived in \cite{Xu} is a little different from \eqref{4.7}. The main difference is that the term involved with $r$ is missing in the former because of those two conditions. The condition number \eqref{4.7} without the term on $r$ and with $L=I_n$ and $\Psi=\beta=\xi=1$ will be equivalent to the one in \cite{Xu}.
}
\end{remark}

Now we set $\mu=\nu=\infty$. In this case, \eqref{4.5} doesn't hold any more. But, considering the properties of the structure matrices $\Phi_{\mathbb{S}_{1}}$ and $\Phi_{\mathbb{S}_{2}}$ (e.g., \cite[Theorem 4.1]{Lijia} and \cite[p. 10]{Rump1}) and the definition of the $\infty$-norm, and noting \eqref{4.400}, we have
\begin{equation}
\label{eq.infty&snorm}
 \left\|\left[ {\begin{array}{*{20}c}
   { {\rm vec}(\frac{\Delta A}{\Psi})}  \\
   { \frac{\Delta b}{\beta}}  \\
\end{array}} \right]\right\|_{\infty}=
    \left\|\left[ {\begin{array}{*{20}c}
   \Phi_{\mathbb{S}_1} &0  \\
   0&\Phi_{\mathbb{S}_2}  \\
\end{array}} \right]\left[ {\begin{array}{*{20}c}
   \frac{\Delta s _1}{\varphi} \\
  \frac{\Delta s_2}{\theta} \\
\end{array}} \right]\right\|_{\infty}=\left\|\left[ {\begin{array}{*{20}c}
   \frac{\Delta s _1}{\varphi} \\
  \frac{\Delta s_2}{\theta} \\
\end{array}} \right]\right\|_{\infty}.
\end{equation}
Substituting \eqref{2.44} and \eqref{eq.infty&snorm} into \eqref{def.scd} and using \eqref{4.200000} yields
\begin{align}
&\kappa^S_{\infty ILS}(A,b)=\mathop {\sup }\limits_{\tiny{\left\|\left[ {\begin{array}{*{20}c}
   \frac{\Delta s _1}{\varphi} \\
  \frac{\Delta s_2}{\theta} \\
\end{array}} \right]\right\|_{\infty}\neq 0}\hfill \atop
\Delta s_1,\varphi \in \mathbb{R}^{k_1}, \Delta s_2, \theta \in \mathbb{R}^{k_2}\hfill}\frac{\left\|{\rm diag}^\ddag(\xi){M_{g' }\left[ {\begin{array}{*{20}c}
   \Phi_{\mathbb{S}_1} &0  \\
   0&\Phi_{\mathbb{S}_2}  \\
\end{array}} \right]\left[ {\begin{array}{*{20}c}
   {\Delta s _1} \\
  {\Delta s_2} \\
\end{array}} \right]}\right\|_{\infty}}{\left\|\left[ {\begin{array}{*{20}c}
   \frac{\Delta s _1}{\varphi} \\
  \frac{\Delta s_2}{\theta} \\
\end{array}} \right]\right\|_{\infty}}\nonumber\\
&\quad\quad\quad\quad\quad =\mathop {\sup }\limits_{\tiny{\left\|\left[ {\begin{array}{*{20}c}
   \frac{\Delta s _1}{\varphi} \\
  \frac{\Delta s_2}{\theta} \\
\end{array}} \right]\right\|_{\infty}\neq 0}\hfill \atop
\Delta s_1,\varphi \in \mathbb{R}^{k_1}, \Delta s_2, \theta \in \mathbb{R}^{k_2}\hfill}\frac{\left\|{\rm diag}^\ddag(\xi){M_{g' }\left[ {\begin{array}{*{20}c}
   \Phi_{\mathbb{S}_1} &0  \\
   0&\Phi_{\mathbb{S}_2}  \\
\end{array}} \right]\mathrm{diag}(\mathrm{vec}(\varphi,\theta))\left[ {\begin{array}{*{20}c}
   \frac{\Delta s _1}{\varphi} \\
  \frac{\Delta s_2}{\theta} \\
\end{array}} \right]}\right\|_{\infty}}{\left\|\left[ {\begin{array}{*{20}c}
   \frac{\Delta s _1}{\varphi} \\
  \frac{\Delta s_2}{\theta} \\
\end{array}} \right]\right\|_{\infty}}.\label{4.33}
\end{align}

From \eqref{4.33}, we have the following theorem.

\begin{theorem}[\bf Structured $\infty$-norm]
Let $A\in \mathbb{S}_1$ and $b \in \mathbb{S}_2$. Then the structured partial condition number under $\infty$-norm for the ILS problem \eqref{1.1} with respect to $L$ is
\begin{eqnarray}\label{eq.sgcd}
\kappa_{\infty ILS}^S (A,b) &=& \left\|{\rm diag}^\ddag(\xi){M_{g' } \left[ {\begin{array}{*{20}c}
   \Phi_{\mathbb{S}_1} &0  \\
   0&\Phi_{\mathbb{S}_2}  \\
\end{array}} \right]\mathrm{diag}(\mathrm{vec}(\varphi,\theta))} \right\|_{\infty}\\
&=&\left\|\frac{{\left|M_{g' } \left[ {\begin{array}{*{20}c}
   \Phi_{\mathbb{S}_1} &0  \\
   0&\Phi_{\mathbb{S}_2}  \\
\end{array}} \right]\right|\left|\begin{bmatrix}
                     \varphi \\
                      \theta \\
                     \end{bmatrix}\right|
}}{{|\xi|}} \right\|_{\infty},
\end{eqnarray}
where $M_{g' }$ is defined by \eqref{2.3}.

In particular, setting $\Psi=A$ and $\beta=b$, i.e., $\varphi=s_1$ and $\theta=s_2$, and $\xi=\left[\|L^Tx(A,b)\|_{\infty},\cdots,\|L^Tx(A,b)\|_{\infty}\right]^T$ or $\xi=L^Tx(A,b)$, we have the structured partial mixed or componentwise condition number
\begin{eqnarray}
\label{eq.smcd}
  \kappa_{mILS}^S (A,b) &=&\frac{\left\|{\left|M_{g' } \left[ {\begin{array}{*{20}c}
   \Phi_{\mathbb{S}_1} &0  \\
   0&\Phi_{\mathbb{S}_2}  \\
\end{array}} \right]\right|\left|\begin{bmatrix}
                     s_1 \\
                      s_2 \\
                     \end{bmatrix}\right|
} \right\|_{\infty}}{{\|L^Tx(A,b)\|_{\infty}}}
\end{eqnarray}
or
\begin{eqnarray}
\label{eq.sccd}
  \kappa_{cILS}^S (A,b) &=&\left\|\frac{{\left|M_{g' } \left[ {\begin{array}{*{20}c}
   \Phi_{\mathbb{S}_1} &0  \\
   0&\Phi_{\mathbb{S}_2}  \\
\end{array}} \right]\right|\left|\begin{bmatrix}
                      s_1 \\
                      s_2 \\
                     \end{bmatrix}\right|
}}{{|L^Tx(A,b)|}} \right\|_{\infty}.
\end{eqnarray}
\end{theorem}

As the special case, the corresponding results for the LLS problem can be obtained.
\begin{corollary}\label{cor100}
Let $A\in \mathbb{S}_1$ and $b \in \mathbb{S}_2$. Then the structured partial condition number under $\infty$-norm for the LLS problem with respect to $L$ is
{\small\begin{align*}
\kappa_{\infty LLS}^S (A,b) = \left\|{{\rm diag}^\ddag(\xi)}{\widetilde{M}_{\widetilde{g}' } \left[ {\begin{array}{*{20}c}
   \Phi_{\mathbb{S}_1} &0  \\
   0&\Phi_{\mathbb{S}_2}  \\
\end{array}} \right]\mathrm{diag}(\mathrm{vec}(\varphi,\theta))} \right\|_{\infty}
=\left\|\frac{{\left|\widetilde{M}_{\widetilde{g}' } \left[ {\begin{array}{*{20}c}
   \Phi_{\mathbb{S}_1} &0  \\
   0&\Phi_{\mathbb{S}_2}  \\
\end{array}} \right]\right|\left|\begin{bmatrix}
                     \varphi \\
                      \theta \\
                     \end{bmatrix}\right|
}}{{|\xi|}} \right\|_{\infty},
\end{align*}}
where $\widetilde{M}_{\widetilde{g}' }$ is defined by \eqref{2.333333}.

In particular, we also have the structured partial mixed and componentwise condition numbers
{\small\begin{align*}
  \kappa_{mLLS}^S(A,b) =\frac{\left\|\left| \left[ {{{\left( {r^T  \otimes (L^T (A^TA)^{ - 1} )} \right)\Pi _{mn}  - x^T  \otimes (L^T A^{\dag})}},{{L^T A^{\dag}}}} \right] \left[ {\begin{array}{*{20}c}
   \Phi_{\mathbb{S}_1} &0  \\
   0&\Phi_{\mathbb{S}_2}  \\
\end{array}} \right]\right|\left|\left[ {\begin{array}{*{20}c}
   s _1 \\
   s_2 \\
\end{array}} \right]\right| \right\|_{\infty}}{\left\|L^Tx(A,b)\right\|_{\infty}} ,\\
   \kappa_{cLLS}^S(A,b)=\left\|\frac{\left| \left[ {{{\left( {r^T  \otimes (L^T (A^TA)^{ - 1} )} \right)\Pi _{mn}  - x^T  \otimes (L^T A^{\dag})}},{{L^T A^{\dag}}}} \right] \left[ {\begin{array}{*{20}c}
   \Phi_{\mathbb{S}_1} &0  \\
   0&\Phi_{\mathbb{S}_2}  \\
\end{array}} \right]\right|\left|\left[ {\begin{array}{*{20}c}
   s _1 \\
   s_2 \\
\end{array}} \right]\right|}{\left|L^Tx(A,b)\right|}\right\|_{\infty}.
\end{align*}}
\end{corollary}

\begin{remark}{\rm
When $L=I_n$, the structured mixed and componentwise condition numbers in Corollary \ref{cor100} will be equivalent to the ones in \cite{Cucker}.}
\end{remark}

\begin{remark}{\rm
All the structures involved in the above results are linear. It is interesting to consider the structured partial condition numbers for the ILS problem under nonlinear structures as done in \cite{Cucker} for the Moore-Penrose inverse and the LLS problem. We will consider this topic in the future research.
}
\end{remark}
\vspace{-6pt}

\section{Connections between ILS and TLS problems}
\label{sec.5}
\vspace{-2pt}
In this section, we consider the partial condition numbers of a special ILS problem, from which we can obtain the corresponding results for the TLS problem.

Let the matrices and vector in the ILS problem \eqref{1.1} be
\begin{eqnarray*}
\widetilde{A}= \left[ {\begin{array}{*{20}{c}}
A\\
B
\end{array}} \right], \ \widetilde{b}= \left[ {\begin{array}{*{20}{c}}
b\\
d
\end{array}} \right],\ \widetilde{J}= \left[ {\begin{array}{*{20}{c}}
I_m&0\\
0&-I_s
\end{array}} \right],
\end{eqnarray*}
where $A\in \mathbb{R}^{m\times n}$ and $b\in \mathbb{R}^{m}$ are the same as the ones in \eqref{1.1}, $B\in \mathbb{R}^{s\times n}$ and $d\in \mathbb{R}^{s}$, and assume that $\widetilde{M}=\widetilde{A}^T\widetilde{J}\widetilde{A}$ is positive definite. Thus, from \eqref{2.1}, it follows that
\begin{align}
&g' (\widetilde{A},\widetilde{b}){\circ}(\Delta \widetilde{A},\Delta \widetilde{b})
=  L^T \widetilde{M}^{ - 1} \left[ {\begin{array}{*{20}{c}}
\Delta A\\
\Delta B
\end{array}} \right]^T \widetilde{J}\widetilde{r} - L^T \widetilde{M}^{ - 1} \widetilde{A}^T \widetilde{J}\left[ {\begin{array}{*{20}{c}}
\Delta A\\
\Delta B
\end{array}} \right]x + L^T \widetilde{M}^{ - 1} \widetilde{A}^T \widetilde{J}\left[ {\begin{array}{*{20}{c}}
\Delta b\\
\Delta d
\end{array}} \right]\nonumber\\
&\quad\quad\quad\quad\quad\quad\quad    =L^T \widetilde{M}^{ - 1} ((\Delta A)^Tr-(\Delta B)^Ts) - L^T \widetilde{M}^{ - 1} (A^T\Delta A-B^T\Delta B)x \nonumber\\
&\quad\quad\quad\quad\quad\quad\quad\quad\  +  L^T \widetilde{M}^{ - 1} (A^T\Delta b-B^T\Delta d).\label{3.2}
\end{align}
In deriving \eqref{3.2}, the result $\widetilde{r}= \left[ {\begin{array}{*{20}{c}}
r\\
s
\end{array}} \right]$ with $r=b-Ax$ and $s=d-Bx$ is used.
As done in \eqref{2.44}, using \eqref{1.6} and \eqref{1.7}, we have
\begin{align}
&g' (\widetilde{A},\widetilde{b}){\circ}(\Delta \widetilde{A},\Delta \widetilde{b})
 =\left[\left(r^T\otimes (L^T \widetilde{M}^{ - 1})\right)\Pi_{mn}-x^T\otimes (L^T \widetilde{M}^{ - 1}A^T), L^T \widetilde{M}^{ - 1}A^T\right]\left[ {\begin{array}{*{20}{c}}
{\rm vec}(\Delta A)\\
\Delta b
\end{array}} \right]\nonumber\\
&\quad\quad+\left[x^T\otimes (L^T \widetilde{M}^{ - 1}B^T)-\left(s^T\otimes (L^T \widetilde{M}^{ - 1})\right)\Pi_{sn} , - L^T \widetilde{M}^{ - 1}B^T\right]\left[ {\begin{array}{*{20}{c}}
{\rm vec}(\Delta B)\\
\Delta d
\end{array}} \right]
.\label{3.3}
\end{align}
Now assume that both $B$ and $d$ are the differentiable functions of $A$ and $b$, i.e., $B$ and $d$ can be written as
\begin{eqnarray*}
B=f_1(A, b),\quad d=f_2(A, b).
\end{eqnarray*}
As a result,
\begin{eqnarray*}
\Delta B=f_1'(A, b)\circ(\Delta A, \Delta b)+\mathcal{O}(\left\|(\Delta A,\Delta b)\right\|_F^2),\quad \Delta d=f_2'(A, b)\circ(\Delta A, \Delta b)+\mathcal{O}(\left\|(\Delta A,\Delta b)\right\|_F^2).
\end{eqnarray*}
Omitting the higher-order terms and using the properties of Kronecker product, the above equations can be written as (e.g.,\cite[p.257]{Horn91})
\begin{eqnarray}\label{3.4}
{\rm vec}(\Delta B)=[M_1,M_2]\left[ {\begin{array}{*{20}{c}}
{\rm vec}(\Delta A)\\
\Delta b
\end{array}} \right],\quad \Delta d=[M_3,M_4]\left[ {\begin{array}{*{20}{c}}
{\rm vec}(\Delta A)\\
\Delta b
\end{array}} \right].
\end{eqnarray}
Substituting \eqref{3.4} into \eqref{3.3} implies
\begin{eqnarray*}
 g' ({A},{b}){\circ}(\Delta {A},\Delta {b}) &=& \left[N_1, N_2\right]\left[ {\begin{array}{*{20}{c}}
{\rm vec}(\Delta A)\\
\Delta b
\end{array}} \right] +\left[N_3 , N_4\right]\left[ {\begin{array}{*{20}{c}}
M_1&M_2\\
M_3&M_4
\end{array}} \right]\left[ {\begin{array}{*{20}{c}}
{\rm vec}(\Delta A)\\
\Delta b
\end{array}} \right]\nonumber\\
&=&\begin{bmatrix}
     N_1+N_3M_1+N_4M_3, & N_2+N_3M_2+N_4M_4 \\
   \end{bmatrix}\left[ {\begin{array}{*{20}{c}}
{\rm vec}(\Delta A)\\
\Delta b
\end{array}} \right],
\end{eqnarray*}
where
\begin{eqnarray}
&&N_1=\left(r^T\otimes (L^T \widetilde{M}^{ - 1})\right)\Pi_{mn}-x^T\otimes (L^T \widetilde{M}^{ - 1}A^T),\quad N_2=L^T \widetilde{M}^{ - 1}A^T,\label{10000}\\
&&N_3=x^T\otimes (L^T \widetilde{M}^{ - 1}B^T)-(s^T\otimes (L^T \widetilde{M}^{ - 1}))\Pi_{sn},\quad N_4=-L^T \widetilde{M}^{ - 1}B^T.\label{10001}
\end{eqnarray}
Thus, analogous to the proof of Theorem \ref{thm.gcils}, we have an expression of the partial unified condition number of this special ILS problem with respect to $L$:
\begin{equation}\label{eq.dftls11}
\kappa_{SILS}(A,b)=\left\|{\rm diag}^\ddag(\xi)\begin{bmatrix}
     N_1+N_3M_1+N_4M_3, & N_2+N_3M_2+N_4M_4 \\
   \end{bmatrix}\mathrm{diag} (\mathrm{vec}(\Psi,\beta)) \right\|_{\mu,\nu}.
\end{equation}


Next, we consider a specific case. That is, set $B=\widetilde{\sigma}_{n+1}I_n$ and $d=0$. Here, $\widetilde{\sigma}_{n+1}$ is the smallest singular value of the matrix $[A,b]$ and is always assumed to be smaller than $\sigma_n$, the smallest singular value of $A$. In this case, $\widetilde{M}=\widetilde{A}^T\widetilde{J}\widetilde{A}=A^TA-\widetilde{\sigma}_{n+1}^2I_n$ is positive definite, and hence the specific ILS problem has the unique solution
\begin{eqnarray*}
x(A,b)=\widetilde{M}^{-1}\widetilde{A}^T\widetilde{J}\widetilde{b}
=(A^TA-\widetilde{\sigma}_{n+1}^2I_n)^{-1}{A}^Tb,
\end{eqnarray*}
which is just the unique solution to the TLS problem expressed as
\begin{equation*}\label{3.7}
{\rm TLS:}\quad\mathop {\min }\limits_{E,\epsilon} \left\|[E,\epsilon]\right\|_F, {\textrm { subject to }} (A+E)x=b+\epsilon.
\end{equation*}
This problem was initially discussed in the seminal paper \cite{Golub} and has many applications \cite{Huf}.

For the above case, from \cite{Bab}, we have
\begin{eqnarray*}
\Delta B=\widetilde{\sigma}_{n+1}^{-1}\frac{r^T(\Delta b-\Delta Ax)}{1+\left\|x\right\|_2^2}I_n+\mathcal{O}(\left\|(\Delta A,\Delta b)\right\|_F^2),
\end{eqnarray*}
and hence
\begin{eqnarray*}
M_1=-\frac{\widetilde{\sigma}_{n+1}^{-1}}{(1+\left\|x\right\|_2^2)}{\rm vec}(I_n)(x^T\otimes r^T),\quad M_2=\frac{\widetilde{\sigma}_{n+1}^{-1}}{(1+\left\|x\right\|_2^2)}{\rm vec}(I_n)r^T.
\end{eqnarray*}
Thus, noting $s=-\widetilde{\sigma}_{n+1}x$, $M_3=0$, $M_4=0$, \eqref{10000}, and \eqref{10001}, and using \eqref{1.6}, \eqref{1.7} and \eqref{1.9}, we obtain
\begin{align*}
 & N_1+N_3M_1=\left(r^T\otimes (L^T \widetilde{M}^{ - 1})\right)\Pi_{mn}-x^T\otimes (L^T \widetilde{M}^{ - 1}A^T)\\
&\quad\quad\quad\quad\quad\quad -\frac{1}{(1+\left\|x\right\|_2^2)}\left(x^T\otimes (L^T \widetilde{M}^{ - 1})+(x^T\otimes (L^T \widetilde{M}^{ - 1}))\Pi_{sn}\right){\rm vec}(I_n)(x^T\otimes r^T),\\
&\quad\quad\quad\quad\ =\left(r^T\otimes (L^T \widetilde{M}^{ - 1})\right)\Pi_{mn}-x^T\otimes (L^T \widetilde{M}^{ - 1}A^T)\\
&\quad\quad\quad\quad\quad\quad-\frac{2L^T \widetilde{M}^{ - 1}x(x^T\otimes r^T)}{(1+\left\|x\right\|_2^2)}\quad \textrm{ by \eqref{1.6} and \eqref{1.7}}\\
&\quad\quad\quad\quad\ =\left(r^T\otimes (L^T \widetilde{M}^{ - 1})\right)\Pi_{mn}-x^T\otimes (L^T \widetilde{M}^{ - 1}A^T)-\frac{2(x^T\otimes (L^T \widetilde{M}^{ - 1}xr^T))}{(1+\left\|x\right\|_2^2)}\ \textrm{ by \eqref{1.9}}
\end{align*}
and
\begin{align*}
&N_2+N_3M_2=L^T \widetilde{M}^{ - 1}A^T+\frac{1}{(1+\left\|x\right\|_2^2)}\left(x^T\otimes (L^T \widetilde{M}^{ - 1})+(x^T\otimes (L^T \widetilde{M}^{ - 1}))\Pi_{sn}\right){\rm vec}(I_n)r^T\\
&\quad\quad\quad\quad\ =L^T \widetilde{M}^{ - 1}A^T+\frac{2L^T \widetilde{M}^{ - 1}xr^T}{(1+\left\|x\right\|_2^2)}.\quad \textrm{ by \eqref{1.6} and \eqref{1.7}}
\end{align*}
Setting $D_\sigma=L^T\widetilde{M}^{-1}\left(A^T+\frac{2xr^T}{1+\left\|x\right\|_2^2}\right)$ and considering \eqref{eq.dftls11} implies an expression of the partial unified condition number for the above specific ILS problem, i.e., the TLS problem with respect to $L$:
\begin{align}\label{TLS}
\kappa_{TLS}(A,b)=\left\|{\rm diag}^\ddag(\xi)\left[\left( \left(r^T\otimes (L^T \widetilde{M}^{ - 1})\right)\Pi_{mn}-x^T\otimes D_\sigma\right),  D_\sigma\right]\mathrm{diag} (\mathrm{vec}(\Psi,\beta)) \right\|_{\mu,\nu}.
\end{align}
Further, if we set $\mu=\nu=2$ and the parameters $\Psi$, $\beta$, and $\xi$ to be positive real numbers, then
\begin{eqnarray}\label{3.8}
\kappa_{2TLS}(A,b)=\frac{\left\|\left[\Psi\left( \left(r^T\otimes (L^T \widetilde{M}^{ - 1})\right)\Pi_{mn}-x^T\otimes D_\sigma\right), \beta D_\sigma\right] \right\|_2}{\xi},
\end{eqnarray}
which reduces to the result in \cite[Proposition 2]{Bab} when $\Psi=\beta=\xi=1$.
\begin{remark}{\rm
Based on \eqref{3.8} or the equivalent expressions, some authors derived the closed formulas and the lower and upper bounds of the normwise condition number for the TLS problem using the singular value decompositions of $A$ and $[A,b]$ and the special properties of the TLS problem (e.g., \cite{Bab,Jia,Lijia}). For the general ILS problem, we have not obtained the corresponding results. 
}
\end{remark}

Now we set $\mu=\nu=\infty$, and $\Psi=A$, $\beta=b$, and $\xi=\left[\|L^Tx(A,b)\|_{\infty},\cdots,\|L^Tx(A,b)\|_{\infty}\right]^T$ or $\xi=L^Tx(A,b)$. Thus, similar to Theorem \ref{thm.m&c}, the partial mixed or componentwise condition number for the TLS problem with respect to $L$ can be obtained
\begin{eqnarray}
\label{eq.mcdtls}
\kappa_{mTLS} (A,b) =\frac{\left\|{\left|\left[ \left(r^T\otimes (L^T \widetilde{M}^{ - 1})\right)\Pi_{mn}-x^T\otimes D_\sigma,  D_\sigma\right]\right|  \left|\mathrm{vec}(A,b)\right| } \right\|_{\infty}}{\|L^Tx(A,b)\|_{\infty}}
\end{eqnarray}
or
\begin{eqnarray}
\label{eq.ccdtls}
\kappa_{cTLS} (A,b) =\left\| \frac{{\left|\left[ \left(r^T\otimes (L^T \widetilde{M}^{ - 1})\right)\Pi_{mn}-x^T\otimes D_\sigma,  D_\sigma\right]\right|  \left|\mathrm{vec}(A,b)\right| }}{|L^Tx(A,b)|} \right\|_{\infty}.
\end{eqnarray}
\begin{remark}{\rm
When $L=I_n$, the mixed and componentwise condition numbers in \eqref{eq.mcdtls} and \eqref{eq.ccdtls} are equivalent to the corresponding ones derived from \cite{Zhou}; see also the discussions in \cite{Lijia,Xie14}.
}

\end{remark}

If the structure of the data in the TLS problem is taken into consideration, as done in Section \ref{sec.4}, we can obtain the structured partial condition numbers for the TLS problem. These results without proof are presented as follows.

The structured partial condition number under 2-norm with the parameters $\Psi$, $\beta$, and $\xi$ being positive real numbers is
\begin{align}
\kappa_{2TLS}^S(A,b)=\frac{\left\|\left[\Psi\left(\left(r^T\otimes (L^T \widetilde{M}^{ - 1})\right)\Pi_{mn}-x^T\otimes D_\sigma\right)\Phi_{\mathbb{S}_1}D_1^{-1},\ {\beta}D_\sigma\Phi_{\mathbb{S}_2} D_2^{-1}\right] \right\|_2}{\xi}.\label{4.10}
\end{align}
Setting $L=I_n$ and $\Psi=\beta=\xi=1$ in \eqref{4.10} leads to the structured normwise condition number for the TLS problem, which is equivalent to the one derived from \cite[(29)]{Lijia}.

The structured partial mixed or componentwise condition number, i.e., the involved norm is $\infty$-norm and the parameters $\Psi=A$, $\beta=b$, and $\xi=\left[\|L^Tx(A,b)\|_{\infty},\cdots,\|L^Tx(A,b)\|_{\infty}\right]^T$ or $\xi=L^Tx(A,b)$ is
\begin{align}
  \kappa_{mTLS}^S(A,b) =\frac{\left\|\left| \left[\left({\left(r^T\otimes (L^T \widetilde{M}^{ - 1})\right)\Pi_{mn}-x^T\otimes D_\sigma}\right)\Phi_{\mathbb{S}_1},\ {D_\sigma}\Phi_{\mathbb{S}_2} \right]\right|\left|\left[ {\begin{array}{*{20}c}
   s _1 \\
   s_2 \\
\end{array}} \right]\right| \right\|_{\infty}}{\left\|L^Tx(A,b)\right\|_{\infty}}  \label{mix_cond}
\end{align}
or
\begin{align}
   \kappa_{cTLS}^S(A,b)=\left\|\frac{\left| \left[\left({\left(r^T\otimes (L^T \widetilde{M}^{ - 1})\right)\Pi_{mn}-x^T\otimes D_\sigma}\right)\Phi_{\mathbb{S}_1},\ {D_\sigma}\Phi_{\mathbb{S}_2} \right] \right|\left|\left[ {\begin{array}{*{20}c}
   s _1 \\
   s_2 \\
\end{array}} \right]\right|}{\left|L^Tx(A,b)\right|}\right\|_{\infty},\label{comp_cond}
\end{align}
where the structured mixed condition number \eqref{mix_cond} with $\Phi_{\mathbb{S}_2}=I_m$, that is, there is no structure requirement on $b$, and $L=I_n$ is equivalent to  \cite[(30) with $\lambda=1$]{Lijia}.

\section{Statistical condition estimates}
\label{sec.6}
In this part, we focus on the methods for estimating the condition numbers under $2$-norm or $\infty$-norm for the ILS problem.
\subsection{Estimating condition number under $2$-norm }
From Theorem \ref{thm.closedform}, we find that the main task to estimate $\kappa_{2ILS}(A,b)$ lies in how to obtain a reliable estimate of the spectral norm of a matrix. This can be carried out by the probabilistic spectral norm estimator \cite{Hochs13}. Meanwhile, an approach based on the SSCE method \cite{Kenney94} can also be applied to estimate $\kappa_{2ILS}(A,b)$ with $L=I_n$. In the following, the brief introductions on these two methods and the corresponding algorithms are presented.

\textbf{PCE method} As mentioned in Section \ref{sec.1}, the probabilistic spectral norm estimator was proposed in \cite{Hochs13}, which provides a reliable estimate of the spectral norm. More precisely, a detailed analysis of the estimator in \cite{Hochs13} showed that the spectral norm of a matrix can be contained in a small interval $[\alpha_1 , \alpha_2]$ with high probability. Here, $\alpha_1$ is the guaranteed lower bound of the spectral norm of the matrix derived by the famous Lanczos bibdiagonalization method \cite{Golub0} and $\alpha_2$ is the probabilistic upper bound of probability at least $1-\varepsilon$ with $\varepsilon\ll 1$ derived by finding the largest zero of a polynomial. Moreover, we can require $\alpha_2/ \alpha_1\leq 1+\delta$ with $\delta$ being a user-chosen parameter.  Based on this estimator, Algorithm \ref{PCE} for estimating the partial condition number \eqref{2.5'} can be devised. 
\begin{algorithm}[htbp]
\caption{Probabilistic condition estimator}\label{PCE}\small
\begin{enumerate}
  \item Compute the matrix
\begin{eqnarray*}
  S=L^TM^{-1}\begin{bmatrix}
                         \Psi\|r\|_2(I_n-\frac{1}{\|r\|^2_2}A^Trx^T), & -\beta A^T, & \Psi \|x\|_2A^T(I_m-\frac{1}{\|r\|_2^2}rr^T) \\
                       \end{bmatrix}
\end{eqnarray*}
and choose a starting random vector $v_0$ from $\mathcal{U}(S_{2m+n-1})$, the uniform distribution over unit sphere $S_{2m+n-1}$ in $R^{2m+n}$.
  \item Compute the guaranteed lower bound $\alpha_1$ and the probabilistic upper bound $\alpha_2$ of $\|S\|_2$ by the probabilistic spectral norm estimator.
  \item Estimate the partial condition number \eqref{2.5'} by
  \begin{eqnarray*}
  \kappa_{p2ILS}(A,b)=\frac{\alpha_1+\alpha_2}{2\xi}.
  \end{eqnarray*}
\end{enumerate}
\end{algorithm}

\begin{remark}{\rm
It is well-known that finding the inverse of a matrix is expensive, and from \cite{Hochs13}, we know that in performing the probabilistic spectral norm estimator, what we really need is the product of the a vector, say $v_0$, with the matrix $S$ or $S^T$, but not the explicit form of $S$. Thus, we can do the following procedure to avoid computing $M^{-1}$. Let
$$D=\begin{bmatrix}
 \Psi\|r\|_2(I_n-\frac{1}{\|r\|^2_2}A^Trx^T), & -\beta A^T, & \Psi \|x\|_2A^T(I_m-\frac{1}{\|r\|_2^2}rr^T) \\
\end{bmatrix}$$
and solve the linear equation $Mx=Dv_0$. Then, let $y=L^Tx$, which is just the product of $v_0$ and $S$, i.e., $Sv_0$. In a similar way, we can compute $S^Tv_0$.

In addition, we can also find that Algorithm \ref{PCE} is applicable to estimating the structured partial condition number \eqref{4.2222} according to the introduction on the probabilistic spectral norm estimator.
}\end{remark}

\textbf{SSCE method} In \cite{Babo}, an approach based on the SSCE method \cite{Gudmun95,Kenney94} is employed to estimate the normwise condition number for the LLS problem, and is showed to perform quite well. Now we apply the approach to estimate the condition number \eqref{2.5} with $L=I_n$. Denote by ${\kappa_{2ILS}}_i(A,b)$ the condition number under 2-norm of the function $z_i^Tx(A,b)$, where $z_i$s are chosen from $\mathcal{U}(S_{n-1})$ and are orthogonal. From \eqref{2.5}, it is seen that
{\small\begin{equation}\label{5.1}
\kappa_{2ILSi} (A,b)  = \frac{ \sqrt{z_i^T M^{ - 1}\left( \Psi^2\left\| r \right\|_2^2 I_n  + \left(\Psi^2\left\| x \right\|_2^2 + \beta ^2 \right)A^T A\right)M^{ - 1}z_i - 2 \Psi^2z_i^T M^{ - 1}xr^T A M^{ - 1}z_i}}{\xi}.
\end{equation}}
The analysis in \cite{Babo} shows that
\begin{eqnarray}\label{5.2}
\kappa_{s2ILS} (A,b) = \frac{\omega_k}{\omega_n}\sqrt{\sum_{i=1}^k {\kappa^2_{2ILS}}_i (A,b)}
\end{eqnarray}
is a good estimate of the condition number \eqref{2.5} with $L=I_n$. In the above expression, $\omega_k$ is the Wallis factor with $\omega_1=1$, $\omega_2={2}/{\pi}$, and
$$\omega_k  =\left\{
  \begin{array}{ll}
    \frac{1\cdot 3 \cdot5 \cdots (k-2)}{2\cdot 4\cdot 6\cdots (k-1)}, & \hbox{for $k$ odd,} \\
    \frac{2}{\pi}\frac{2\cdot 4\cdot 6\cdots (k-2)}{3\cdot 5\cdot 7\cdots (k-1)}, & \hbox{for $k$ even,}
  \end{array}
\right.\textrm{ when } k>2.$$
It can be approximated by
\begin{eqnarray}\label{5.3}
\omega_k\approx \sqrt{\frac{2}{\pi(k-\frac{1}{2})}}
\end{eqnarray}
with high accuracy. In summary, we can propose Algorithm \ref{SSCE}.

\begin{algorithm}[htbp]{
\caption{An approach based on small-sample statistical condition estimation}\label{SSCE}\small
\begin{enumerate}
  \item Generate $k$ vectors $z_1,\cdots, z_k$ from $\mathcal{U}(S_{n-1})$, and orthonormalize these vectors using the QR facotization.
  \item For $i=1,\cdots, k$, compute ${\kappa_{2ILS}}_i (A,b)$ by \eqref{5.1}.
  \item Approximate $\omega_k$ and $\omega_n$ by \eqref{5.3} and estimate the condition number \eqref{2.5} with $L=I_n$ by \eqref{5.2}.
\end{enumerate}}
\end{algorithm}

\subsection{Estimating condition numbers under $\infty$-norm}
For the partial mixed and componentwise condition numbers, we consider the SSCE method \cite{Gudmun95,Kenney94}, which has been used to estimate the condition numbers for the linear systems, the LLS problem,  the matrix equations et al.  (e.g., \cite{Babo,R_Diao12,R_Kenney98,R_KenneyLS98, Laub08}). As done in the aforementioned references, we devise Algorithm \ref{m&cSSCE} to estimate the condition numbers $\kappa_{mILS}(A,b)$ and $\kappa_{cILS}$ in \eqref{eq.mcdils} and \eqref{eq.ccdils}.\\
\begin{algorithm}[htbp]
\caption{Small-sample statistical condition estimation}\label{m&cSSCE}\small
\begin{enumerate}
  \item Let $t=m(n+1)$. Generate $k$ vectors ${{z}}_1,\cdots, {{z}}_k$ from $\mathcal{U}(S_{t-1})$, and orthonormalize these vectors using the QR facotization.
  \item Compute $u_i=M_{g'}{z}_i$, and estimate the  partial mixed and componentwise condition numbers in \eqref{eq.mcdils} and \eqref{eq.ccdils} by
  \begin{equation*}
  \kappa_{smILS}(A,b)=\frac{\|\kappa_{ssce}\|_{\infty}}{\|L^Tx(A,b)\|_{\infty}}, \quad \kappa_{scILS}(A,b)=\left\|\frac{\kappa_{ssce}}{L^Tx(A,b)}\right\|_{\infty},
  \end{equation*} where $\kappa_{ssce}=\frac{\omega_k}{\omega_t}\left|\sum^{k}_{i=1}|u_i|^2\right|^{\frac{1}{2}}$, and the power and square root operation are performed on each entry of $u_i$, $i=1,\cdots, k$.
\end{enumerate}
\end{algorithm}
\vspace{-6pt}

\section{Numerical experiments}
\label{sec.7}
\vspace{-2pt}
Three numerical examples are presented in this section. The first two are used to illustrate the reliability of the statistical condition estimators presented by Algorithms \ref{PCE} and \ref{SSCE}, and Algorithm \ref{m&cSSCE}, respectively, and the third one is used to compare the structured partial condition numbers and the unstructured ones. In these examples, for simplicity, we set $\Psi=\beta=\xi=1$ and the matrix $L$ to be the identity matrix.

\begin{example}\label{examp1}
{\rm
In this example, the ILS problem \eqref{1.1} is generated as follows. First, form the matrix $A\in\mathbb{R}^{m\times n}$ by
\begin{equation*}
    A=\left[ {\begin{array}{*{20}{c}}
U_p & 0\\
0 & U_q
\end{array}} \right]\begin{bmatrix}
         D \\
         0 \\
       \end{bmatrix}
V,\; U_p=I_p-2u_pu_p^T,\;  U_q=I_q-2u_qu_q^T,\;\mathrm{and} \; V=I_n-2vv^T,
\end{equation*}
where $u_p\in\mathbb{R}^{p}$, $u_q\in\mathbb{R}^{q}$ and $v\in\mathbb{R}^{n}$ are unit random vectors and $D=n^{-l}\mathrm{diag}(n^{l}, (n-1)^{l},\cdots,1^l)$. It is easy to find that the condition number of $A$, i.e., $\kappa(A)=\left\|A\right\|_2\left\|A^\dag\right\|_2$, is $n^l$. Then, set the solution $x$ to be $x=(1, 2^2,\cdots,n^2)$ and $b=Ax+r$ with $r$ being random vector of 2-norm $\rho$, i.e., $\|r\|_2=\rho$.

In practical implementation, we set $m=200$, $n=120$, and $p=140$. It can be easily checked that $A^TJA$ is positive definite for this setting. For Algorithm \ref{PCE}, we choose the parameters $\delta= 0.01$ and $\epsilon=0.001$. In this case, the inequalities $\alpha_1\leq \left\|S\right\|_2\leq \alpha_2$ hold with a probability at least $99.9\%$, and $\alpha_1$ and $\alpha_2$ satisfy the inequality  $\alpha_2/\alpha_1\leq 1.01$. For Algorithm \ref{SSCE}, we set $k=3$. By varying the condition number of $A$ and the 2-norm of the residual vector $r$, we use $500$ ILS problems for each pair of $\kappa(A)$ and $\rho$ to test the performance of the aforementioned two algorithms. To show the efficiency of statistical condition estimators clearly, we define the ratios between the estimate and the exact value as follows
\begin{eqnarray*}
  r_p &=& \kappa_{p2ILS}(A,b)/\kappa_{2ILS}(A,b),\;r_s=\kappa_{s2ILS}(A,b)/\kappa_{2ILS}(A,b),
\end{eqnarray*}
and report the mean and variance of these ratios in Table \ref{Table1}.

\begin{table}[!htb]
\centering
\caption{The efficiency of statistical condition estimates in Algorithms \ref{PCE} and \ref{SSCE}}\label{Table1}
\footnotesize
\begin{tabular}{c|c|cccccccc}
  \hline
    &$\kappa(A)$ & \multicolumn{2}{c}{$n^{0}$} & \multicolumn{2}{c}{$n^{3}$}  \\
\hline
         $\rho$    &  & mean & variance  &  mean  & variance   \\
\hline
 $10^{-4}$ &$r_p$&  1.000e+00&  6.845e-11&  1.000e+00&  6.671e-11\\
           &$r_s$&  1.197e+01&  4.149e-19&  1.023e+00&  1.764e-01\\
\hline
$10^{-2}$ &$r_p$&   1.000e+00&  8.104e-11 & 1.000e+00&  5.585e-11 \\
          &$r_s$&   1.197e+01&  4.233e-15 & 1.034e+00&  1.796e-01 \\
\hline
$10^{0}$ &$r_p$&    1.000e+00&  8.346e-11 & 1.000e+00&  8.690e-11 \\
         &$r_s$&    1.197e+01&  4.346e-11 & 9.723e-01&  1.618e-01\\
\hline
$10^{2}$ &$r_p$&    1.001e+00&  7.953e-11&  1.000e+00&  8.530e-11 \\
         &$r_s$&    1.197e+01&  3.990e-07 & 1.032e+00&  1.801e-01 \\
\hline
$10^{4}$ &$r_p$&    1.000e+00&  1.057e-10 & 1.000e+00&  8.682e-11\\
         &$r_s$&    1.138e+01&  3.071e-03&  1.025e+00&  1.743e-01 \\
\hline
&$\kappa(A)$  & \multicolumn{2}{c}{$n^{6}$} & \multicolumn{2}{c}{$n^{9}$} \\

\hline
 $10^{-4}$ &$r_p$&  1.001e+00&  2.310e-06&  1.000e+00&  3.566e-08\\
           &$r_s$&  1.253e+00&  1.313e-01&  1.442e+00&  1.197e-01\\
\hline
$10^{-2}$ &$r_p$ & 1.002e+00 & 1.055e-06 & 1.000e+00 & 3.248e-07\\
          &$r_s$ & 1.188e+00 & 1.385e-01&  1.354e+00&  1.443e-01\\
\hline
$10^{0}$ &$r_p$ & 1.000e+00 & 1.371e-11 & 1.000e+00 & 2.729e-08\\
         &$r_s$&  1.079e+00 & 1.480e-01 & 1.174e+00 & 1.290e-01\\
\hline
$10^{2}$ &$r_p$ & 1.000e+00 & 1.319e-11&  1.000e+00 & 2.987e-08\\
         &$r_s$ & 1.084e+00 & 1.662e-01&  1.146e+00  & 1.509e-01\\
\hline
$10^{4}$ &$r_p$&  1.000e+00 & 1.298e-11 & 1.000e+00 & 3.505e-08\\
         &$r_s$ & 1.034e+00&  1.531e-01&  1.158e+00&  1.419e-01\\
\hline
\end{tabular}

\end{table}

According to the explanation in \cite[Chapter 15]{Hig02}, `an estimate of the condition number that is correct to within a factor 10 is usually acceptable, because it is the magnitude of an error bound that is of interest, not its precise value', the results in Table \ref{Table1} show that both Algorithms \ref{PCE} and \ref{SSCE} can give reliable estimates of the condition number under 2-norm in most cases. In comparison, Algorithm \ref{PCE} performs better and more stable, but Algorithm \ref{SSCE} may behave bad when $\kappa(A)=1$. This latter phenomenon also appears in estimating the normwise condition number of the LLS problem; see \cite{Babo} for an explanation.

}\end{example}

\begin{example}\label{examp2}
{\rm This example is constructed according to \cite{Boa} and Example \ref{examp1}. That is, the matrix $A$ is formed as
 \begin{equation*}
 A=\left[ {\begin{array}{*{20}{c}}
Q_1DU\\
\frac{1}{2}Q_2DU
\end{array}} \right],
 \end{equation*}
where $Q_1\in {\mathbb{R}^{p \times n}}$, $Q_2\in {\mathbb{R}^{q \times n}}$, and $U\in {\mathbb{R}^{n \times n}}$ are the random orthogonal matrices and $D\in {\mathbb{R}^{n \times n}}$ is a diagonal matrix with diagonal elements distributed exponentially from $\kappa^{-1}$ to 1. Then, set the solution $x$ and the residual vector $r$ as done in Example \ref{examp1}.

In the numerical experiments, we set $m = 120$, $n = 50$, $p = 70$, and $k=3$. For each pair of $\kappa$ and $\rho$, $200$ ILS problems are generated to test the performance of Algorithm \ref{m&cSSCE}. The numerical results on mean and variance of the ratios between the statistical condition estimate and
the exact condition number defined by
\begin{eqnarray*}
  r_m &=& \kappa_{smILS}(A,b)/\kappa_{mILS}(A,b)\textrm{ and }r_c=\kappa_{scILS}(A,b)/\kappa_{cILS}(A,b)
\end{eqnarray*}
are reported in Tables \ref{Table2}, which suggest that Algorithm \ref{m&cSSCE} is effective and
reliable in estimating the mixed and componentwise condition numbers.

\begin{table}[!htb]
\centering
\caption{The efficiency of statistical condition estimates in Algorithm \ref{m&cSSCE}}\label{Table2}
\scriptsize
\begin{tabular}{c|c|llllllll}
  \hline
   & $\kappa(A)$ & \multicolumn{2}{c}{$10^{2}$} & \multicolumn{2}{c}{$10^{6}$}  \\
\hline
          $\rho$   &  & mean & variance  &  mean  & variance   \\
\hline
 $10^{-4}$ &$r_m$&  1.024e+00&  3.625e-02&  1.409e+00&  1.483e-01\\
           &$r_c$&  6.754e-01&  9.795e-02&  1.092e+00 & 2.135e-01\\
\hline
$10^{-2}$ &$r_m$&  1.008e+00&  3.280e-02 & 1.258e+00&  1.562e-01\\
          &$r_c$&  6.957e-01&  8.891e-02 & 1.066e+00 & 2.045e-01\\
\hline
$10^{0}$ &$r_m$&   1.056e+00 & 4.725e-02 & 1.349e+00 & 2.046e-01\\
         &$r_c$&   6.366e-01 & 8.669e-02 &  1.219e+00 & 2.614e-01\\
\hline
$10^{2}$ &$r_m$&  9.272e-01 & 3.973e-02 & 1.389e+00 & 2.614e-01\\
         &$r_c$&  6.096e-01 & 6.096e-02&  1.225e+00 & 2.910e-01\\
\hline
$10^{4}$ &$r_m$&   1.137e+00 &  8.303e-02&  1.470e+00 & 2.715e-01\\
         &$r_c$&   7.891e-01 & 1.157e-01 & 1.234e+00  & 3.070e-01\\
\hline
   & $\kappa(A)$  & \multicolumn{2}{c}{$10^{10}$} & \multicolumn{2}{c}{$10^{12}$} \\
\hline
 $10^{-4}$ &$r_m$&    1.542e+00&  2.580e-01&  1.616e+00&  3.315e-01\\
           &$r_c$&    1.332e+00&  3.758e-01&  1.412e+00&  3.355e-01\\
\hline
$10^{-2}$ &$r_m$&  1.581e+00 & 4.304e-01&  1.620e+00&  4.649e-01\\
          &$r_c$&   1.465e+00 & 3.770e-01 & 1.569e+00 & 4.732e-01\\
\hline
$10^{0}$ &$r_m$&    1.669e+00 & 3.949e-01 & 1.726e+00 & 4.543e-01\\
         &$r_c$&    1.589e+00 & 3.714e-01 & 1.629e+00 & 4.647e-01\\
\hline
$10^{2}$ &$r_m$&  1.646e+00 & 4.038e-01&  1.733e+00  &  5.482e-01\\
         &$r_c$&   1.586e+00 & 4.050e-01&  1.706e+00 & 6.246e-01\\
\hline
$10^{4}$ &$r_m$&    1.627e+00&  4.258e-01&  1.727e+00&  4.311e-01\\
         &$r_c$&    1.593e+00 & 5.324e-01 & 1.622e+00 & 4.712e-01\\
\hline
\end{tabular}
\end{table}
}\end{example}

\begin{example}\label{examp3}
{\rm The matrix $A$ in this example is formed as
$ A=\left[B^T,\
\frac{1}{2}B^T \right]^T$,
where $B\in {\mathbb{R}^{n \times n}}$ is a nonsymmetric gaussian random Toeplitz matrix generated by the Matlab function ${\bf toeplitz}(c,r)$ with $c={\bf randn}(n,1)$ and $r={\bf randn}(n,1)$, and the solution $x$ and the residual vector $r$ are the same as the ones in Example \ref{examp1}. For the above setting, $m=2n$ and the structure on $b$ is not considered. Meanwhile, we set $p=q=n$ in $J$. In this case, $A^TJA$ is always positive definite when $B$ is nonsingular.

In the practical experiments, we set $n=60$, and generate $200$ ILS problems for each $\rho$. The numerical results on the ratios defined by
\begin{eqnarray*}
r_N = \kappa_{2ILS}(A,b)/\kappa_{2ILS}^S(A,b),  r_M = \kappa_{mILS}(A,b)/\kappa_{mILS}^S(A,b), r_C=\kappa_{cILS}(A,b)/\kappa_{cILS}^S(A,b)
\end{eqnarray*}
are presented in Table \ref{Table3}. 
These results confirm the analysis in Remark \ref{Remark}. Also, we can find that in some cases the unstructured condition number under 2-norm is much larger than the structured one. 

\begin{table}[!htb]
\centering
\caption{Comparisons of the structured condition numbers and the unstructured ones}\label{Table3}
\footnotesize
\begin{tabular}{c|c|ccccc}
  \hline
  & \backslashbox{ratios}{$\rho$}  & $10^{-4}$ &$10^{-2}$ & $10^{0}$ & $10^{2}$ &$10^{4}$ \\
\hline
            &$r_N$&  8.8414&  8.6376&  8.1939&  8.9524&  8.2121\\
mean        &$r_M$&  4.0977&  4.0108&  3.9511 & 4.1248&  5.5935\\
           &$r_C$&  4.3096&  4.2435&  4.3426 & 4.3549&  5.4572\\
\hline
            &$r_N$&  45.7575&   32.6965&   25.0709&   34.6592&  17.5130\\
max        &$r_M$&  8.8546&    8.0559&    7.6172&    8.1457&  10.4139\\
           &$r_C$&  8.5969&    8.3016&    8.7231&    8.0207&  10.4133\\
\hline
\end{tabular}
\end{table}


}\end{example}
\vspace{-6pt}

\section{Conclusion}
\vspace{-2pt}
This paper first presents an explicit expression of the partial unified condition number of the ILS problem. 
Then, the explicit expressions of the partial normwise, mixed and componentwise condition numbers are obtained. These results generalize the corresponding ones for the LLS problem in \cite{Ari,Bab09} and improve the corresponding ones for the ILS problem in \cite{Boa, Li14}. Corresponding to the unstructured partial condition numbers, the structured ones are also derived, which generalize and improve the corresponding ones for the LLS problem in \cite{Cucker,Xu}. Furthermore, we consider the condition numbers for the TLS problem from the view of the  ILS problem and recover and generalize some results given in \cite{Bab,Lijia}.  As far as we know,  it is the first time to investigate the condition numbers for the TLS problem in this way.
Finally, the statistical estimates of the derived condition numbers and the corresponding algorithms are provided. Numerical experiments show that these estimates are efficient and reliable. Meanwhile, a numerical example also confirms that the structured condition numbers are indeed tighter than the unstructured ones.
\vspace{-6pt}

\section*{Acknowledgments}
{
The authors would like to thank Dr. Michiel E. Hochstenbach for providing Matlab program of the probabilistic spectral norm estimator.}

\end{document}